\documentclass[letterpaper, 10 pt, conference]{ieeeconf}
\IEEEoverridecommandlockouts
\usepackage{graphicx}
\usepackage{float}
\usepackage{array}
\usepackage{amsmath,amssymb}
\usepackage{eepic}
\usepackage{epsfig}
\usepackage{subfigure}
\usepackage{calc}
\usepackage{amssymb}
\usepackage{amstext}
\usepackage{amsmath}
\usepackage{multicol}
\usepackage{pslatex}
\usepackage{float}
\usepackage{listings}
\usepackage{cite}
\newtheorem{example}{Example}

\newtheorem{theorem}{Theorem}

\newtheorem{assumption}{Assumption}
\DeclareMathOperator{\sign}{sign}

\title{\LARGE \bf
A short survey on nonlinear models of the classic Costas loop: \\
rigorous derivation and limitations of the classic analysis.
\thanks{
        Accepted to American Control Conference (ACC) 2015 (Chicago, USA)
}
}

\author{Best R.E.$^{1}$, Kuznetsov N.V.$^{2,3}$, Kuznetsova O.A.$^{3}$, Leonov G.A.$^{3}$,
Yuldashev M.V.$^{3}$, Yuldashev R.V.$^{3}$
\thanks{$^{1}$
        Best Engineering company, Oberwil, Switzerland
        }%
\thanks{$^{2}$
        Mathematical Information Technology Dept., University of Jyv\"{a}skyl\"{a}, Jyv\"{a}skyl\"{a}, Finland
        {\footnotesize (corresponding author: nkuznetsov239@gmail.com)}
        }%
\thanks{$^{3}$
        Faculty of Mathematics and Mechanics, Saint-Petersburg State University, Saint-Petersburg, Russia
        }%
}

\usepackage[update,prepend]{epstopdf}
\begin{document}

\maketitle
\thispagestyle{empty}
\pagestyle{empty}

%%%%%%%%%%%%%%%%%%%%%%%%%%%%%%%%%%%%%%%%%%%%%%%%%%%%%%%%%%%%%%%%%%%%%%%%%%%%%%%%
\begin{abstract}

Rigorous nonlinear analysis of the physical model of Costas loop ---
a classic phase-locked loop (PLL) based circuit for carrier recovery,
is a challenging task.
Thus for its analysis, simplified mathematical models and numerical simulation are widely used.
In this work a short survey on nonlinear models of the BPSK Costas loop,
used for pre-design and post-design analysis, is presented.
Their rigorous derivation and limitations of classic analysis are discussed.
It is shown that
the use of simplified mathematical models,
and the application of non rigorous methods of analysis (e.g., simulation and linearization)
may lead to wrong conclusions concerning
the performance of the Costas loop physical model.

\end{abstract}

%%%%%%%%%%%%%%%%%%%%%%%%%%%%%%%%%%%%%%%%%%%%%%%%%%%%%%%%%%%%%%%%%%%%%%%%%%%%%%%%
\section{INTRODUCTION}

The Costas loop is a classic phase-locked loop (PLL) based circuit
for carrier recovery \cite{Costas-1956,Costas-1962-patent}.
In this paper the classic analog Costas loop \cite{Costas-1956,Best-2007},
used for Binary Phase Shift Keying signals (BPSK) is considered.
Costas loop is essentially a nonlinear control system
and its {\it physical model} is described by
a nonlinear non-autonomous system of discontinuous differential equations
(\emph{mathematical model in the signal space}), whose rigorous analytical is a difficult task.
Thus, in practice, numerical simulation, simplified mathematical models, and linear analysis
are widely used for the analysis of PLL based circuits
(see, e.g., \cite{KiharaOE-2002,Bianchi-2005-book,Best-2007,PedersonM-2008-book,Tranter-2010-book,Talbot-2012-book}).

In the following it is shown that
1) the use of simplified mathematical models,
and
2) the application of non rigorous methods of analysis (e.g., a simulation)
may lead to wrong conclusions about the performance of the Costas loop \emph{physical model}.

To demonstrate this, the operation of the Costas loop will be considered in details.

\section{Classical engineering consideration of the Costas loop operation}

\begin{figure}[thpb]
\centering
 \includegraphics[width=0.4\textwidth]{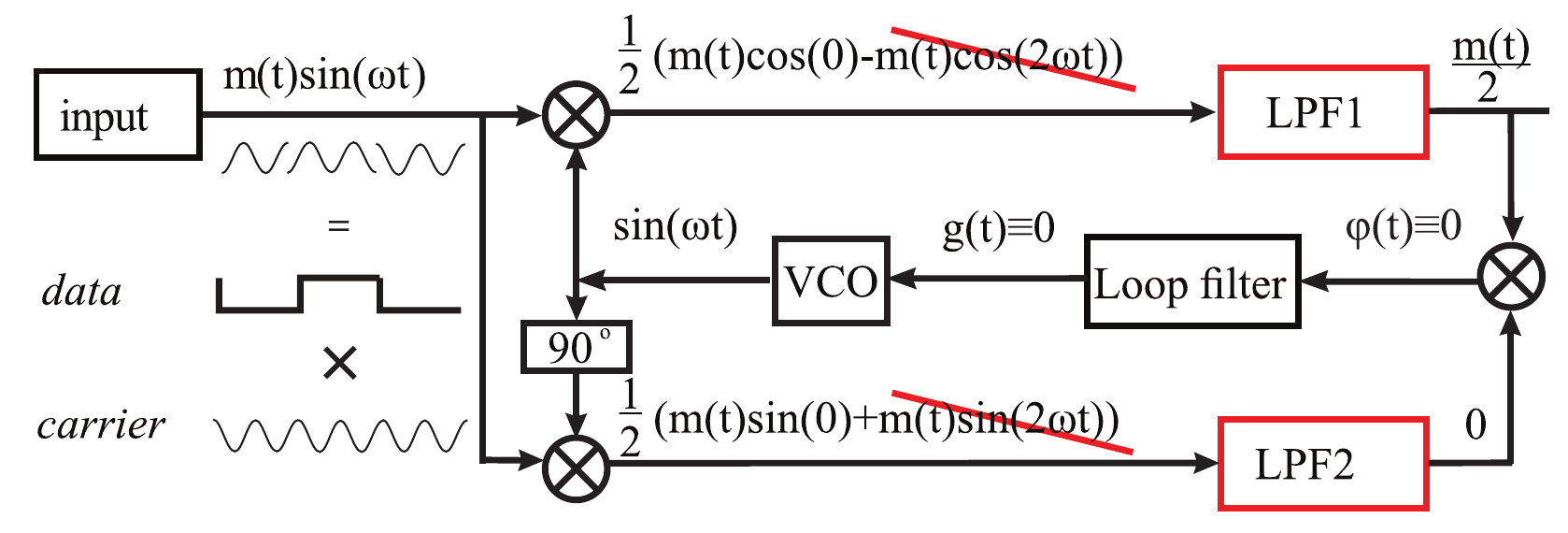}
\caption{Costas loop is locked
(the case of equal phases of input carrier and free running VCO output):
there is no phase difference.
}
\label{costas_locked}
\vspace{-0.4cm}
\end{figure}

The operation of the Costas loop is considered first in the locked state (see Fig.~\ref{costas_locked}),
hence the frequency of the carrier is identical with the frequency
of the VCO (Voltage-Controlled Oscillator);
further it is assumed that both of these signals are sinusoidal.

The input signal is a BPSK signal,
which is the product of a transferred binary data ($m(t) \in\{ \pm 1 \}$ for any $t$)
and the harmonic carrier $\sin(\omega t)$ with a high frequency $\omega$.
Since the Costas loop is considered to be locked,
the VCO output signal is synchronized with the carrier
(i.e. there is no phase difference between these two signals).
The input signal is multiplied (multiplier block ($\otimes$))
by the VCO signal on the upper branch and by the VCO signal, shifted by $90^{\circ}$, on the lower branch.
Therefore on the multipliers' outputs one has
\(
 \varphi_1(t) =
 \frac{1}{2}
 \big(
 m(t) - m(t)\cos(2\omega t)
 \big),
 \varphi_2(t) = \frac{1}{2}
 \big(
 m(t)\sin(2\omega t)
 \big).
\)

Consider the low-pass filters (LPF1 and LPF2) operation.
\begin{assumption}\label{as-twice-frequency}
  Signals components, whose frequency is about twice the carrier frequency,
  do not affect the synchronization of the loop
  (since they are completely suppressed by the low-pass filters).
\end{assumption}
\begin{assumption}\label{as-lpf-initstate}
  Initial states of the low-pass filters do not affect the synchronization of the loop
  (since for the properly designed filters, the impact of filter's initial state on its output
  decays exponentially with time).
\end{assumption}
\begin{assumption}\label{as-lpf-data}
  The data signal $m(t)$ does not affect the synchronization of the loop.
\end{assumption}

Assumptions~\ref{as-twice-frequency},\ref{as-lpf-initstate}, and \ref{as-lpf-data}
together lead to the concept of so-called \emph{ideal low-pass filter}, which
completely eliminates all frequencies above the cutoff frequency (Assumption~\ref{as-twice-frequency})
while passing those below unchanged (Assumptions~\ref{as-lpf-initstate},\ref{as-lpf-data}).
In the classic engineering theory of the Costas loop it is assumed that
the low-pass filters LPF1 and LPF2 are ideal low-pass filters.

Since in Fig.~\ref{costas_locked} the loop is in lock,
i.e. the transient process is over and the synchronization is achieved,
by Assumptions~\ref{as-twice-frequency},\ref{as-lpf-initstate}, and \ref{as-lpf-data}
for the outputs $g_{1,2}(t)$ of the low-pass filters LPF1 and LPF2 one has
\(
   g_1(t) = \frac{1}{2} m(t), \ g_2(t) = 0.
\)
Thus, the upper branch works as a demodulator and the lower branch works as a phase-locked loop.

Since after a transient process there is no phase difference,
a control signal at the input of VCO,
which is used for VCO frequency adjustment to the frequency of input carrier signal,
has to be zero: $g(t) = 0$.
In the general case when the carrier frequency $\omega$ and
a free-running frequency $\omega_{\text{free}}$ of the VCO are different,
after a transient processes
the control signal at the input of VCO has to be non-zero constant:
$g(t) = const$,
and a constant phase difference $\theta_\Delta$ may remain.
\smallskip
\smallskip
\smallskip

\begin{figure}[thpb]
\centering
 \includegraphics[width=0.4\textwidth]{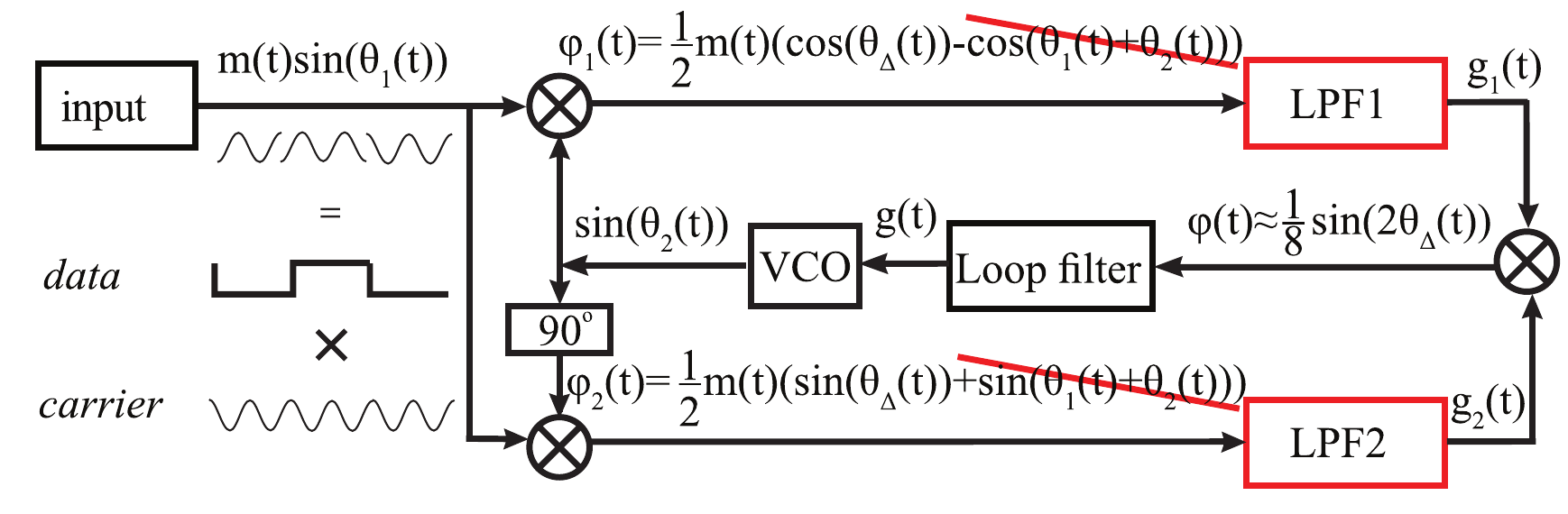}
 \caption{Costas loop is out of lock: there is time-varying phase difference.
 }
\label{costas_out}
\end{figure}

Consider the Costas loop before synchronization
(see Fig.~\ref{costas_out}).
Here the phase difference
$
  \theta_{\Delta}(t) = \theta_1(t)-\theta_2(t)
$
varies over time, because the loop has not yet acquired lock
(frequencies or phases of the carrier and VCO are different).

In this case, using Assumption~\ref{as-twice-frequency},
the signals $\varphi_{1,2}(t)$ can be approximated as
\begin{equation}\label{phi12-approx}
\begin{aligned}
  & \varphi_1(t)
  %= \frac{1}{2}m(t)
  %\big(\cos(\theta_1(t)-\theta_2(t)) - \cos(\theta_1(t)+\theta_2(t)) \big)\\
  %&
  \approx \frac{1}{2}m(t)\cos(\theta_{\Delta}(t)),
  %\\
  & \varphi_2(t)
  %=  \frac{1}{2}m(t)
  %\big(\sin(\theta_1(t)-\theta_2(t)) + \sin(\theta_1(t)+\theta_2(t))  \big)\\
  %&
  \approx \frac{1}{2}m(t)\sin(\theta_{\Delta}(t)).
\end{aligned}
\end{equation}
Approximations \eqref{phi12-approx} depend on the phase difference of signals,
i.e. two multiplier blocks ($\otimes$) on the upper and lower branches
operate as phase detectors.
The obtained expressions \eqref{phi12-approx} with $m(t)\equiv1$
coincide with well-known (see, e.g., \cite{Viterbi-1966,Best-2007})
phase detector characteristic of the classic PLL with multiplier/mixer phase-detector
for sinusoidal signals.

By Assumptions~\ref{as-lpf-initstate} and \ref{as-lpf-data}
the low-pass filters outputs can be approximated as
\begin{equation}\label{g12-approx}
\begin{aligned}
  & g_1(t)
  \approx \frac{1}{2}m(t)\cos(\theta_{\Delta}(t)),
  & g_2(t)
  \approx \frac{1}{2}m(t)\sin(\theta_{\Delta}(t)).
\end{aligned}
\end{equation}
For $m^2(t) \equiv 1$, the input of the loop filter is
\begin{equation}\label{loop-filter-input-approx}
 \varphi(t) = g_1(t)g_2(t) \approx \varphi(\theta_{\Delta}(t)) = \frac{1}{8} \sin(2\theta_{\Delta}(t)).
\end{equation}
Such an approximation
is called \emph{phase detector characteristic of the Costas loop}.

Since an ideal low-pass filter is hardly realized,
its use in the mathematical analysis
requires additional justification.

Thus, the impact of the low-pass filters on the lock acquisition process must be studied rigorously.
Also, the following caveats should be acknowledged.

{\bf Caveat to Assumption~\ref{as-twice-frequency}}.
While Assumption~\ref{as-twice-frequency}
is reasonable from a practical point of view,
its use in the analysis of Costas loop requires further consideration
(see, e.g., \cite{PiqueiraM-2003,KuznetsovKLNYY-2015-ISCAS}).
Various averaging methods
(see, e.g., \cite{MitropolskyB-1961})
allow one to justify Assumption~\ref{as-twice-frequency}
and obtain conditions under which it can be used rigorously
(see \cite{LeonovKYY-2012-TCASII,LeonovKYY-2015-SIGPRO}).

{\bf Caveat to Assumption~\ref{as-lpf-initstate}}.
Since in Fig.~\ref{costas_out} the loop is out of lock,
i.e. synchronization is not achieved,
low-pass filters' initial states cannot be ignored
and must be taken into account.
Therefore for rigorous consideration of low-pass filters influence
one has to use a rigorous mathematical model of the low-pass filters
instead of approximations \eqref{g12-approx}.
Since low-pass filter LPF2 is mostly used
to indicate synchronization status
and low-pass filter LPF1 is mostly used for data demodulation,
the effect of nonzero initial state of filter
on transient processes will be discussed below for the loop filter,
which is used to provide synchronization.

\smallskip

{\bf Caveat to Assumption~\ref{as-lpf-data}}.
The low-pass filters can not operate perfectly at
the moments of changing $m(t)$,
therefore the data pulse shapes
are no longer ideal rectangular pulses after filtration.
%due to distortion, created by the low-pass filters.
One known related effect is called false-locking:
while for $m(t) \equiv const$
the loop acquires lock and proper synchronization of the carrier and VCO frequencies,
for time-varying $m(t) \neq const$
the loop may seem to acquire lock without proper synchronization of the frequencies (false lock)
\cite{Olson-1975,Stensby-1989}.  %,Simon-1978,Lindsey-1978
To avoid such undesirable situation one may
try to choose loop parameters in such a way that
the synchronization time is less than the time between changes in the data signal $m(t)$
or to modify the loop (see, e.g., \cite{Olson-1975}).
%Another way is to perform the nonlinear nonlocal analysis of the loop
%(see, e.g., \cite{Stensby-1989}) %,Stensby-2002
%to identify unsuitable parameters.
\smallskip

The relation between the input $\varphi(t)$
and the output $g(t)$ of the loop filter has the form %\cite{Thede-2005}
\begin{equation}\label{loop-filter}
 \begin{aligned}
 & \dot x = A x + b \varphi(t),
 \ g(t) = c^*x + h\varphi(t),
 \end{aligned}
\end{equation}
where $A$ is a constant matrix,
the vector $x(t)$ is the loop filter state,
$b,c$ are constant vectors, h is a number.
The control signal $g(t)$ is used to adjust the VCO frequency to
the frequency of the input carrier signal
\begin{equation} \label{vco first}
   \dot\theta_2(t) = \omega_2(t) = \omega_2^{\text{free}} + Lg(t).
\end{equation}
Here $\omega_2^{free}$ is the free-running frequency of the VCO
and $L$ is the VCO gain.
The solution of \eqref{loop-filter} with initial data $x(0)$
(the loop filter output for the initial state $x(0)$) is as follows
\begin{equation}\label{loop-filter-int}
 \begin{array}{c}
 g(t,x(0)) = \alpha_0(t,x(0)) +
 \int\limits_0^t
 \gamma(t - \tau)\varphi(\tau)
 {\rm d}\tau
 + h \varphi(t),
 \end{array}
\end{equation}
where $\gamma(t - \tau)=c^*e^{A(t-\tau)}b + h$
is the impulse response of the loop filter
and $\alpha_0(t,x(0))= c^*e^{At}x(0)$
is the zero input response of the loop filter,
i.e. when the input of the loop filter is zero.

\begin{assumption}[analog of Assumption~\ref{as-lpf-initstate}]\label{as-lf-initstate}
  Zero input of loop filter $\alpha_0(t,x(0))$
  does not affect the synchronization of the loop
  (one of the reasons is that
  $\alpha_0(t,x(0))$ is an exponentially damped function for a stable matrix $A$).
\end{assumption}
\smallskip

Consider a constant frequency of the input carrier:
\begin{equation}\label{omega1-const}
   \dot\theta_1(t) = \omega_1(t) \equiv \omega_1,
\end{equation}
and introduce notation: $\omega_{\Delta}^{\text{free}} = \omega_1-\omega_2^{free}$.
Then Assumption~\ref{as-lf-initstate}
allows one to obtain the classic mathematical model of PLL-based circuit
in signal's phase space
(see the classic Viterbi's book \cite{Viterbi-1966}):
\begin{equation}\label{mathmodel-class-simple}
 \begin{aligned}
   & \dot\theta_{\Delta} =
   \omega_{\Delta}^{\text{free}}-L\int_0^t
   \gamma(t - \tau)\varphi(\theta_{\Delta}(\tau)){\rm d}\tau
   -Lh \varphi(\theta_{\Delta}(t)).
 \end{aligned}
\end{equation}
Since $\varphi(\theta_{\Delta})$ from \eqref{loop-filter-input-approx}
is an odd function and has period $\pi$, one has
$(\omega_{\Delta}^{\text{free}},\theta_{\Delta}(t))
\rightarrow \big(-\omega_{\Delta}^{\text{free}},-\theta_{\Delta}(t))$
and $\theta_{\Delta}(t) \rightarrow \theta_{\Delta}(t)+\pi k$
do not change \eqref{mathmodel-class-simple}.
Therefore in $\eqref{mathmodel-class-simple}$ one can consider
only nonnegative $\omega_{\Delta}^{\text{free}}$
($|\omega_{\Delta}^{\text{free}}|$ called \emph{frequency deviation})
and $\theta_{\Delta}(0) \in [-\pi/2,\pi/2)$.
The classic engineering task is to find the following sets:
the hold-in range includes such $|\omega_{\Delta}^{\text{free}}|$
that \eqref{mathmodel-class-simple}
has a stationary state $\theta_{\Delta}(t) \equiv \theta^{eq}_{\Delta}$
which is locally stable
(local stability, i.e. for some $\theta_{\Delta}(0)$);
the pull-in range includes such $|\omega_{\Delta}^{\text{free}}|$
that any solution of \eqref{mathmodel-class-simple}
is attracted to one of the stationary state $\theta^{eq}_{\Delta}$
(global stability, i.e. for all $\theta_{\Delta}(0)$).

%\newpage
{\bf Caveat to Assumption~\ref{as-lf-initstate}}.
While Assumption~\ref{as-lf-initstate} allows
one to introduce one dimensional stability ranges, defined only by $|\omega_{\Delta}^{\text{free}}|$,
$\alpha_0(t,x(0))$ may affect the synchronization of the loop and stability ranges.
For rigorous study
one has to consider multi-dimensional stability domains,
taking into account $(\omega_{\Delta}^{\text{free}},x(0),x_1(0),x_2(0))$,
and explain their relations with the classic engineering ranges \cite{LeonovK-2014}
(e.g., it is of utmost importance for cycle slips study \cite{AscheidM-1982,ErshovaL-1983}
and lock-in range definition).

Note that while Assumption~\ref{as-lpf-initstate}
is explained at the beginning Viterbi's classical book \cite{Viterbi-1966}
for the stable matrices $A$ only,
further in the book, various filters with marginally stable
matrices are also studied
(e.g. filter -- perfect integrator where $A=0$).

It is also interesting that for model $\eqref{mathmodel-class-simple}$ with $h=0$
(see, e.g., eq.~(2.18) in \cite{Viterbi-1966})
the initial difference between frequencies $|\dot \theta_{\Delta}(0)| = |\omega_{\Delta}(0)|$ is equal
to the frequency deviation $|\omega_{\Delta}^{\text{free}}|$.
So these terms are often used in place of each other,
which is not correct for $x(0) \neq 0$, $h \neq 0$ or
a non odd function $\varphi(\theta_{\Delta})$.
\smallskip

While the classic model \eqref{mathmodel-class-simple}
may be useful at the stage of post-design analysis
(when the input and the VCO output are considered only
and the parameters are known only approximately),
for the pre-design analysis
(when all the parameters of the loop can be chosen precisely) one may use
more informative models considered in the next sections.

\section{Nonlinear models of Costas loop}
The relation between the inputs $\varphi_{1,2}(t)$ and the outputs $g_{1,2}(t)$
of linear low-pass filters is as follows %\cite{Thede-2005}
\begin{equation}\label{LPF}
 \begin{aligned}
 & \dot x_{1,2} = A_{1,2} x_{1,2} + b_{1,2}\varphi_{1,2}(t),
 \ g_{1,2}(t) = c_{1,2}^*x_{1,2}.
 \end{aligned}
\end{equation}
Here $A_{1,2}$ are constant stable (all eigenvalues have negative real part) matrices,
the vectors $x_{1,2}(t)$ are low-pass filters' states,
$b_{1,2}$ and $c_{1,2}$ are constant vectors,
the vectors $x_{1,2}(0)$ are initial states of the low-pass filters.
For the loop filter one can consider more the general equation \eqref{loop-filter}
with a proportional term.

Taking into account \eqref{LPF}, \eqref{loop-filter}, and \eqref{vco first},
one obtains the \emph{mathematical model in the signal space},
describing the \emph{physical model} of BPSK Costas loop,:
\begin{equation}
\label{prev diff eq}
 \begin{aligned}
 & \dot x_1 = A_1 x_1 + b_1 m(t)\sin(\theta_1(t))\sin(\theta_2), \\
 & \dot x_2 = A_2 x_2 + b_2 m(t)\sin(\theta_1(t))\cos(\theta_2), \\
 & \dot x = A x + b (c_1^* x_1) (c_2^* x_2), \\
 & \dot \theta_2 = \omega_2^{free} + L(c^*x) + Lh(c_1^* x_1) (c_2^* x_2).
 \end{aligned}
\end{equation}
Note that Assumptions 1-4 are not used in the derivation of system \eqref{prev diff eq}.

The \emph{mathematical model in the signal space} \eqref{prev diff eq}
is a nonlinear nonautonomous discontinuous differential system,
so in general case its analytical study is a difficult task
even for the continuous case when $m(t)\equiv const$.
Besides it is a slow-fast system, so its numerical study
(corresponding to the SPICE level simulation)
is rather complicated for the high-frequency signals.
The problem is that it is necessary
to consider simultaneously both very fast time scale
of the signals $\sin(\theta_{1,2}(t))$
and slow time scale of phase difference
$\theta_{\Delta}(t)$,
therefore a very small simulation time-step
must be taken over a very long total simulation period
\cite{Goyal-2006}. %Abramovitch-2008-1, Abramovitch-2008-2

To overcome these problems,
in place of using Assumption~2
one can apply averaging methods
\cite{MitropolskyB-1961,LeonovKYY-2015-SIGPRO} %Samoilenko-2004-averiging,
and consider \emph{a simplified mathematical model in signal's phase space}.
However, this requires the consideration of a constant data signal (Assumption~\ref{as-lpf-data})
and constant frequency of input carrier \eqref{omega1-const}, i.e.
\[
  m(t)\equiv 1, \quad \theta_1(t) = \omega_1 t + \theta_1(0).
\]
In this case \eqref{prev diff eq} is equivalent to
\begin{equation}
 \label{diff_eq_0}
 \begin{aligned}
 & \dot{x_1} = A_1 x_1
 + b_1 \sin(\omega_1 t + \theta_1(0))
 \sin(\omega_1 t + \theta_1(0) + \theta_{\Delta}), \\
 & \dot{x_2} = A_2 x_2
 + b_2 \sin(\omega_1 t + \theta_1(0))
 \cos(\omega_1 t + \theta_1(0) + \theta_{\Delta}), \\
 & \dot{x} = A x + b (c_1^* x_1) (c_2^* x_2), \\
 & \dot\theta_{\Delta} = \omega_{\Delta}^{\text{free}} - L(c^*x) - Lh(c_1^* x_1) (c_2^* x_2).\\
 \end{aligned}
\end{equation}

Here the initial (at $t=0$)  difference of the frequencies has the form
\begin{equation} \label{freq-init-ss}
   \dot \theta_{\Delta}(0)=\omega_{\Delta}(0)=\omega_{\Delta}^{\text{free}} - Lc^*x(0) - Lh c_1^* x_1(0) c_2^* x_2(0).
\end{equation}

Assuming that the input carrier is a high-frequency signal
(i.e. $\omega_1 > \omega_{min}$ is sufficiently large), one can
consider small parameter
$   \varepsilon = \frac{1}{\omega_1}$.
 Denote
     $ \tau = \omega_1 t$.
 Then system \eqref{diff_eq_0} can be represented in the following way
\begin{equation} \label{standard form}
   \frac{dz}{d\tau} = \varepsilon F(z,\tau),  \quad z = (x_1,x_2,x,\theta_{\Delta})^*.
 \end{equation}
 %In the classic averaging theory
 %such a form of system is called a standard form.
 Consider the averaged equation \eqref{standard form}
 \begin{equation} \label{general averaged}
 \frac{dy}{d\tau} = \varepsilon \bar F(y), \ \bar F(y) = \frac{1}{2\pi}\int_0^{2\pi} F(y,\tau) d\tau.
 \end{equation}

 Suppose, $D$ is a bounded domain,
 containing the point $z_0 = (x_1(0),x_2(0),x(0),\theta_{\Delta}(0))$.
 Consider solutions
 $z(\tau,\varepsilon)$ and $y(\tau,\varepsilon)$
 with the initial data $z_0=y_0$.
 In this case there exists a constant $T$ such that $z(\tau,\varepsilon)$ and $y(\tau,\varepsilon)$
 remain in the domain $D$ for $0 \leq \tau \leq \frac{T}{\varepsilon}$.
 Define $\varepsilon_{max} = \frac{1}{\omega^{min}}$.

  \begin{theorem} \cite{MitropolskyB-1961}
    \label{thm firs order periodic aver}
    Suppose that $\varepsilon_{max}$, $D$, and $T$ are as above.
     Then there exists a constant $c > 0$ such that
     $
       ||y(t,\varepsilon) - z(t,\varepsilon)|| < c\varepsilon
     $
     for $0 \leq \varepsilon \leq \varepsilon_{max}$ and $0 \leq t \leq \frac{T}{\varepsilon}$.
  \end{theorem}
  %Proof of this theorem can be found in \cite{Sanders-2007} (p. 31).

 {\bf Remark}. Time $T$ has to be defined by the time of the transient processes.
 The theorem does not suggest how to define $T$ by the loop parameters,
 so $T$ is supposed to be estimated experimentally.
 Also note that the averaged system can be considered
 only if
 \begin{equation}\label{close-freq}
   \begin{aligned}
   & 0< \omega^{min} \leq \omega_1, \omega_{2}(t), \\
   & |\omega_1 - \omega_2(t)| \leq \omega^{max}_{\Delta} \quad \forall t \in [0,T],
   \end{aligned}
 \end{equation}
 where $\omega^{max}_{\Delta}$ does not depend on $\omega^{min}$
 and $\omega^{min}$ is sufficiently large.
\bigskip

\begin{figure}[thpb]
\centering
 \includegraphics[width=0.4\textwidth]{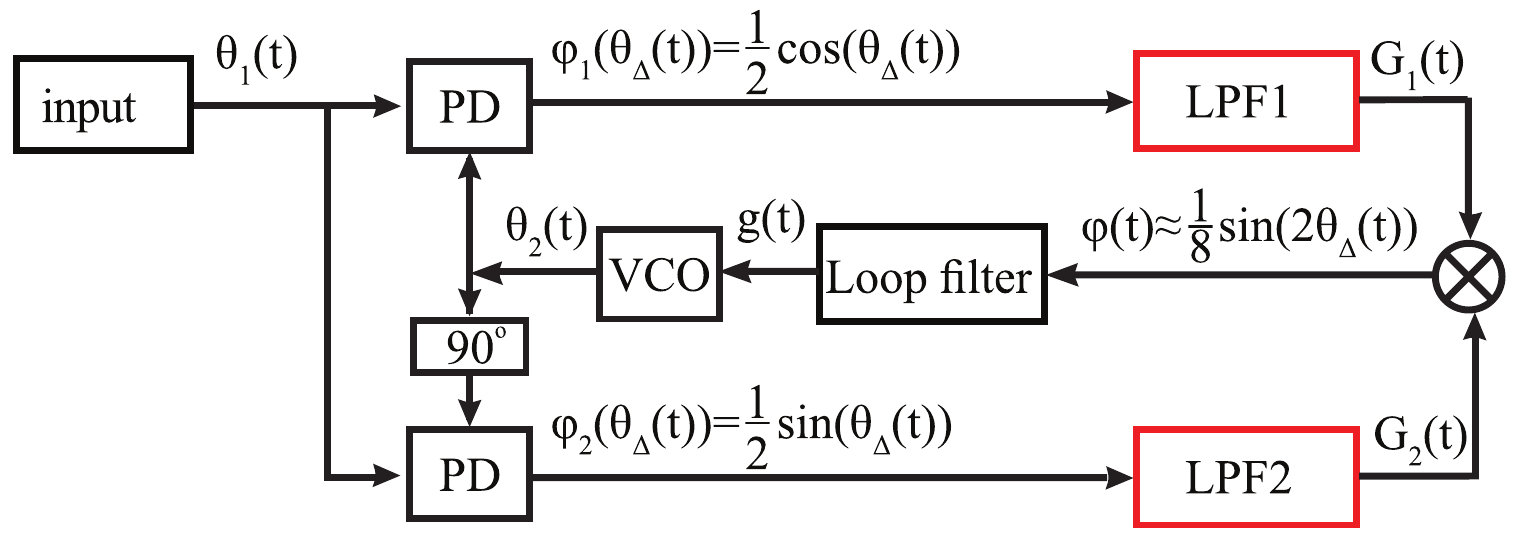}
\caption{Mathematical model of BPSK Costas loop in signal's phase space}
\label{costas_phase_1}
\end{figure}
\vspace{-0.3cm}

In this case the averaged system \eqref{general averaged}
gives a \emph{mathematical model of BPSK Costas loop in signal's phase space} (see Fig.~\ref{costas_phase_1}):
\begin{equation}\label{averaged}
 \begin{aligned}
 & \dot{x_1} = A_1 x_1 + \frac{b_1}{2}\cos(\theta_{\Delta}), \quad
  \dot{x_2} = A_2 x_2 + \frac{b_2}{2}\sin(\theta_{\Delta}), \\
 & \dot{x} = A x + b (c_1^* x_1) (c_2^* x_2), \\
 & \dot\theta_{\Delta} = \omega_{\Delta}^{\text{free}} - L(c^*x) - Lh(c_1^* x_1) (c_2^* x_2).
 \end{aligned}
\end{equation}
Here, under condition \eqref{close-freq}, Theorem~\ref{thm firs order periodic aver}
provides the closeness
solutions of averaged \eqref{averaged} and original \eqref{prev diff eq} systems.
Remark that in \eqref{averaged} the initial difference between frequencies $\omega_{\Delta}(0)$
coincides with \eqref{freq-init-ss}.

Mathematical model \eqref{averaged} does not coincide with that,
considered in the classic works
\cite{Costas-1962-patent,Viterbi-1966}.
The classical mathematical model of the Costas loop can be obtained from this model
under additional assumptions.

\subsection{Classical model of BPSK Costas loop}

Having solved the first two equations of system \eqref{averaged}, one obtains

\(%begin{equation}
 %\label{lpf12-full}
 \begin{aligned}
 & c_1^* x_1 = \alpha_1(t,x_1(0)) +
 \int_{0}^{t}\gamma_1(t - \tau)
 \frac{1}{2}\cos(\theta_{\Delta}(\tau))
 d\tau, \\
 & c_2^* x_2 = \alpha_2(t,x_2(0)) +
 \int_{0}^{t}\gamma_2(t - \tau)
 \frac{1}{2}\sin(\theta_{\Delta}(\tau))
 d\tau. \\
 \end{aligned}
\)%end{equation}

Here the impact of filter's initial state on its output decays exponentially with time
(see Assumption~\ref{as-lpf-initstate}) and one can consider
 $t>t_0$, such that $\alpha_{1,2}(t,x_{1,2}(0)) = O(\frac{1}{\omega^{min}})$.

Since the low-pass filters LPF1 and LPF2 are assumed to be ideal
(see Assumptions~\ref{as-twice-frequency}-\ref{as-lpf-initstate}),
under conditions \eqref{close-freq}
their operation can be formalized as
\[ %begin{equation} \label{lpf_conditions}
\begin{aligned}
 &
 \int_{t_0}^t \gamma_{1,2}(t-\tau)\sin(\theta_{\Delta}(\tau))d\tau
  =
 \sin\big(\theta_{\Delta}(t)\big) + O(\frac{1}{\omega^{min}}).
 %,\  \dot \theta_{\Delta}(t) < \omega^{max}_{\Delta},
\end{aligned}
\] %end{equation}
Then %from \eqref{lpf12-full}
for $t>t_0$ one has
\begin{equation}
 \label{zero initials systems}
 \begin{aligned}
 & c_1^* x_1(t) =
 %\alpha_1(t) +
 %\int\limits_{0}^{t}\gamma_1(t - \tau)
 %\frac{1}{2}\cos(\theta_{\Delta}(\tau)) d\tau  = \\
 %& =
 \dfrac{1}{2}\cos\big(\theta_{\Delta}(t)\big) + O(\frac{1}{\omega^{min}}), \\
 & c_2^* x_2(t) =
 %= \alpha_2(t) +
 %\int\limits_{0}^{t}\gamma_2(t - \tau)
 %\frac{1}{2}\sin(\theta_{\Delta}(\tau)) d\tau = \\
 %& =
 \dfrac{1}{2}\sin\big(\theta_{\Delta}(t)\big) + O(\frac{1}{\omega^{min}}). \\
 \end{aligned}
\end{equation}
Applying \eqref{zero initials systems} to \eqref{averaged}, for $t>t_0$ one obtains
\begin{equation}
 \label{sys-with-O}
 \begin{aligned}
 &\dot{x}  = A x + b \frac{1}{8}\sin(2\theta_{\Delta}) + O(\frac{1}{\omega^{min}}), \\
 & \dot\theta_{\Delta} = \omega_{\Delta}^{\text{free}}- Lc^*x- Lh \frac{1}{8}\sin(2\theta_{\Delta})
 - O(\frac{1}{\omega^{min}}).\\
 \end{aligned}
\end{equation}

\begin{assumption}[Corollary of Assumptions 1-3)]\label{as-averaged}
Solutions of system \eqref{averaged}
are close to solutions of the following system
\begin{equation}
 \label{final_system}
 \begin{aligned}
 & \dot{x} = A x + b \varphi(\theta_{\Delta}), \\
 & \dot\theta_{\Delta} = \omega_{\Delta}^{\text{free}} - Lc^*x - Lh\varphi(\theta_{\Delta}),
 \end{aligned}
\end{equation}
i.e. $O(\frac{1}{\omega^{min}})$ in system \eqref{averaged} can be neglected.
\end{assumption}

\begin{figure}[thpb]
\centering
 \includegraphics[width=0.35\textwidth]{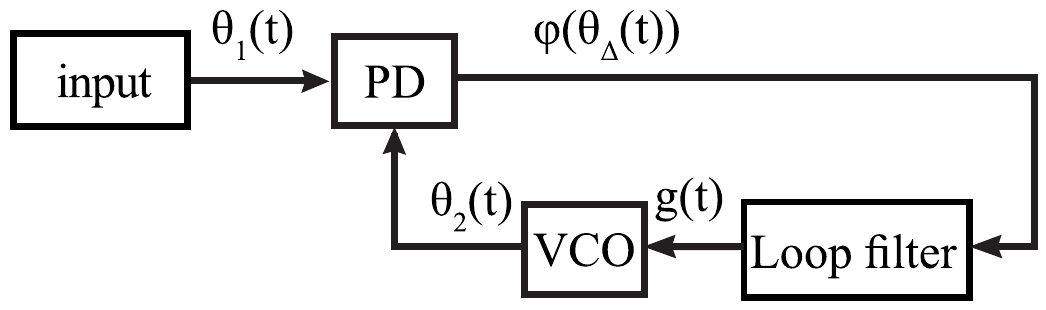}
\caption{Classical (simplified) mathematical model of BPSK Costas loop in signal's phase space}
\label{costas_phase_1_short}
\end{figure}

System \eqref{final_system} corresponds
to the block-diagram shown in Fig.~\ref{costas_phase_1_short},
where $\varphi(\theta_{\Delta})$ is
the phase detector characteristic of the Costas loop for sinusoidal signals.
Note that here the phase detector operation includes
the operations of three multipliers, phase shift element $90^{\circ}$, LPF1, and LPF2.
System \eqref{final_system} with $h=0$ and $x(0)=0$
corresponds to system \eqref{mathmodel-class-simple}.

{\bf Caveat to Assumption~\ref{as-averaged}}.
For rigorous justification of Assumption~\ref{as-averaged}
one may analyze stability conditions for system \eqref{final_system}
(see, e.g., criteria of stability in the large for the pendulum-like systems \cite{LeonovK-2014}).
\smallskip

\begin{figure}[H]
  \centering
  \includegraphics[width=0.4\textwidth]{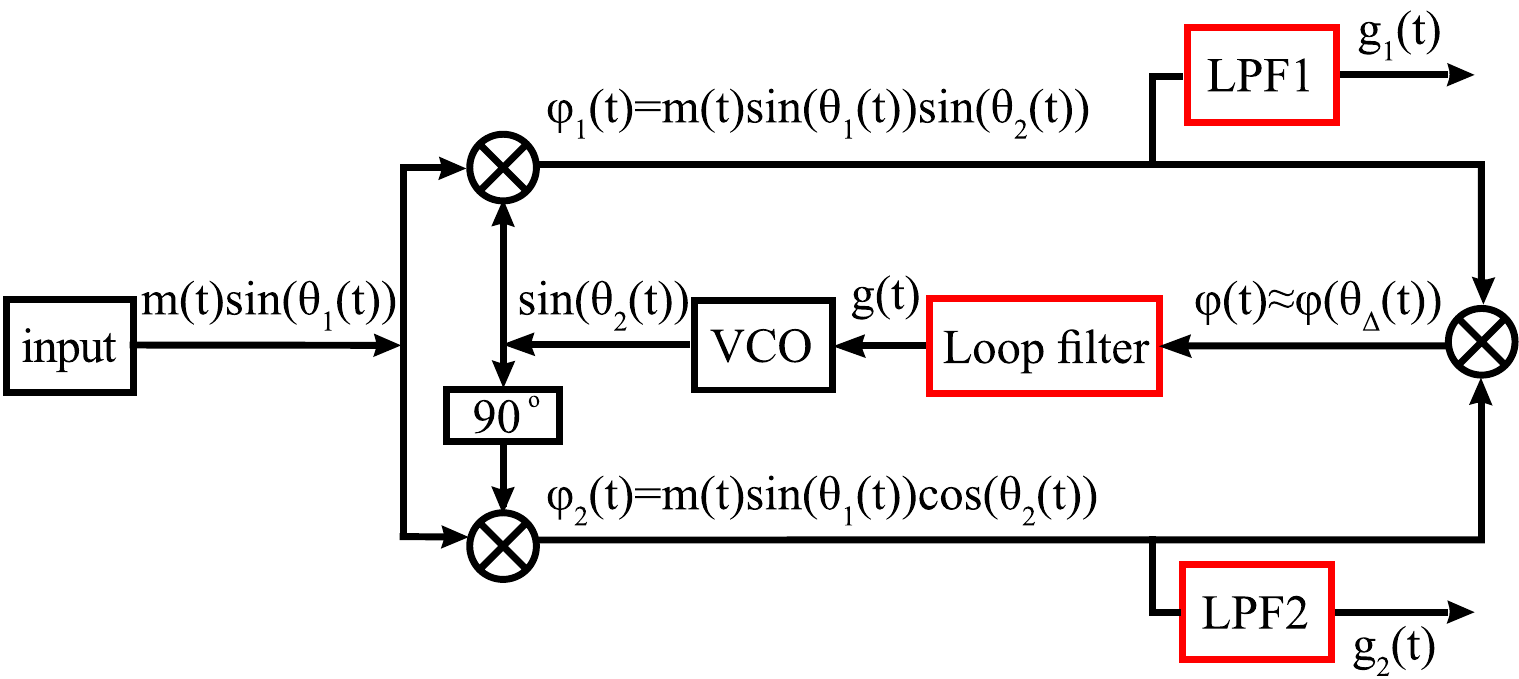}
  \caption{A modification of BPSK Costas loop}
  \label{simplfied_costas}
\end{figure}

In some applications
a modification of the BPSK Costas loop, shown in Fig.~\ref{simplfied_costas}, is used.
Here low-pass filters LPF1 and LPF2 do not affect
the operation of the modified BPSK Costas loop (i.e. Assumption~\ref{as-lpf-initstate} is not needed);
the input of the loop filter does not depend on the data signal $m(t)$
(i.e. Assumption~\ref{as-lpf-data} is not needed):

$\varphi(t)=\varphi_1(t)\varphi_2(t) = \sin(\theta_1(t))\sin(\theta_2(t))\sin(\theta_1(t))\cos(\theta_2(t))$.

Taking into account \eqref{loop-filter} and \eqref{vco first},
one obtains a \emph{mathematical model in the signal space}
describing the \emph{physical model} of modified BPSK Costas loop:
\begin{equation} \label{two_phase_equations}
 \begin{aligned}
 & \dot{x} = A x + b\varphi(t), \quad
  \dot\theta_{\Delta} = \omega_{\Delta}^{\text{free}} -  Lc^*x - Lh\varphi(t).
 \end{aligned}
\end{equation}

The averaged system \eqref{two_phase_equations}
gives a \emph{mathematical model of modified BPSK Costas loop in signal's phase space}
(see Fig.~\ref{costas_phase_1_short}) which corresponds to \eqref{final_system}.
Here, under condition \eqref{close-freq}, from Theorem~\ref{thm firs order periodic aver}
it follows that the solutions of systems \eqref{two_phase_equations} and \eqref{final_system}
are close.
The approach suggested in
\cite{KuznetsovLYY-2011-IEEE,LeonovKYY-2012-TCASII,ourPatent-2013-fin,LeonovKYY-2015-SIGPRO}
allows one to find an approximation $\varphi(\theta_{\Delta})$
in the case of non-constant frequency of the carrier
or non sinusoidal signals.

\section{Summary of the models}
\begin{figure}
\centering
 \includegraphics[width=0.4\textwidth]{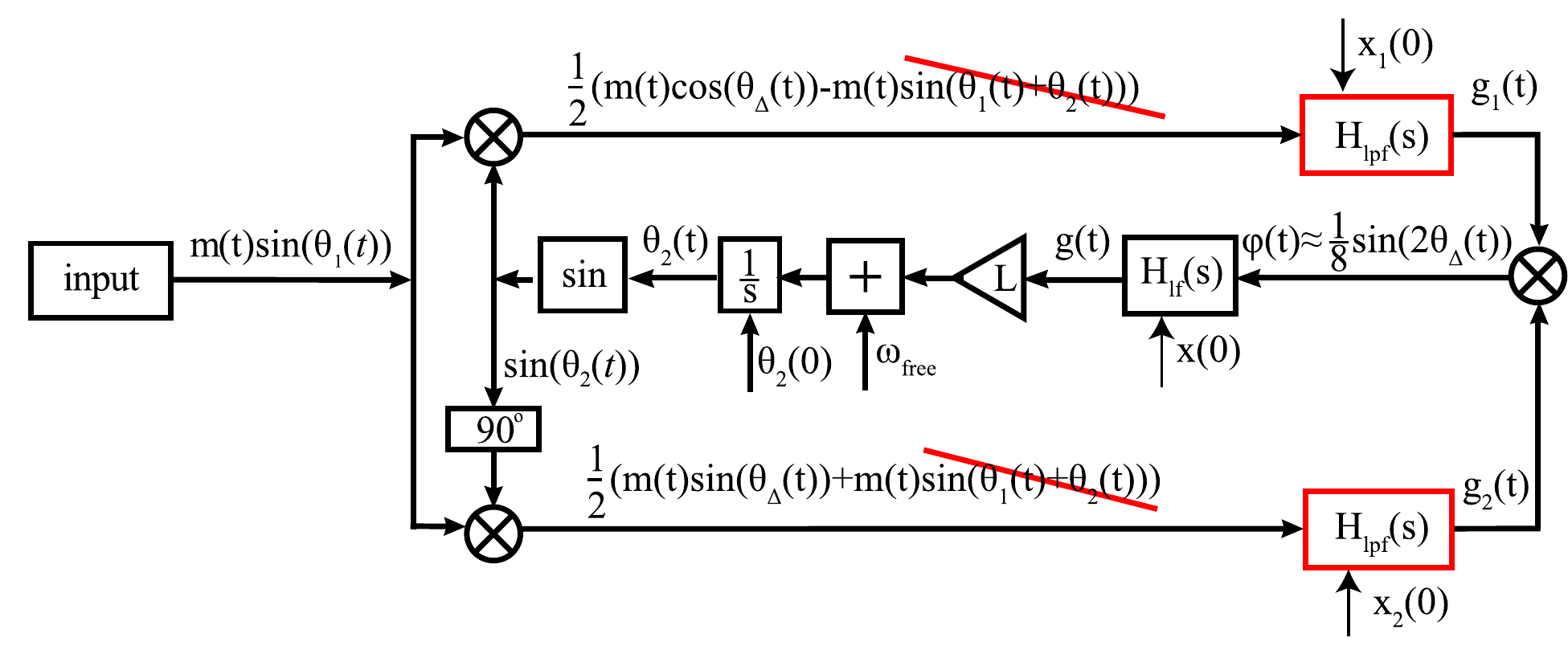}
 \caption{Block-diagram of Costas loop mathematical model in the signal space
 described by transfer functions and initial conditions}
\label{physical_sim_model}
\end{figure}

\begin{itemize}
  \item \emph{physical model} (with data signal and low-pass filters)
        or its \emph{mathematical model in the signal space}
        (Fig.~\ref{costas_out} or Fig.~\ref{physical_sim_model} and system \eqref{prev diff eq});
  \item \emph{simplified mathematical model in the signal space}
       (without data signal and with low-pass filters))
       (system \eqref{diff_eq_0} and Fig.~\ref{costas_out} or Fig.~\ref{physical_sim_model} with $m(t) \equiv 1$);
  \item \emph{modified physical model} (with data signal and with external low-pass filters)
        or its \emph{mathematical model in the signal space}
       (Fig.~\ref{simplfied_costas} and system \eqref{two_phase_equations});
\end{itemize}
%and models in signal's phase space
\begin{itemize}
  \item
   \emph{simplified mathematical model in signal's phase space} (without data signal and with low-pass filters)
   (system \eqref{averaged} and Fig.~\ref{costas_phase_1});
  \item
   \emph{classic mathematical model in signal's phase space}
   (without data signal and without low-pass filters)
   (system \eqref{final_system} and Fig.~\ref{costas_phase_1_short}).
\end{itemize}

\section{Costas loop simulation and counterexamples to the Assumptions}
Note once more that various simplifications
and analysis of linearized models of control systems
may result in incorrect conclusions.
At the same time the attempts to justify analytically the reliability of conclusions,
based on engineering approaches, and rigorous study of nonlinear models
are quite rare (see, e.g., tutorial \cite{Abramovitch-2002}).

Further it is demonstrated that the use of
the above engineering Assumptions
requires further study and rigorous justification.
Next examples shows that for the same parameters
the behaviors of considered models
may be substantially different from one another.

\subsection{Simulation examples}
Next the following parameters are used in simulation:
low-pass filters transfer functions $H_{lpf}(s) = \frac{2}{s/\omega_3+1}$, $\omega_3 = 1.2566*10^6$
and corresponding parameters in system \eqref{loop-filter}  are
$A_{1,2} = -\omega_3$, $b_{1,2} = 1$, $c_{1,2} = \omega_3$;
loop filter transfer function $H_{lf}(s) = \frac{\tau_2 s + 1}{\tau_1 s}$,
$\tau_2 = 3.9789*10^{-6}$, $\tau_1 = 2*10^{-5}$,
and corresponding parameters in system \eqref{loop-filter}  are
$A = 0$, $b = 1$, $c = \frac{1}{\tau_1}$, $h = \frac{\tau_2}{\tau_1}$;
carrier frequency $\omega_1=2*\pi*400000$;
VCO input gain $L=4.8*10^6$;
and carrier initial phase $\theta_2(0)=\theta_1(0)=0$.

\begin{example}[double frequency and averaging]
In Fig.~\ref{averaged_non_averaged}
it is shown that Assumption~\ref{as-twice-frequency} may not be valid
if conditions for the application of Theorem~\ref{thm firs order periodic aver}
are violated:
mathematical model signal's phase space (see Fig.~\ref{costas_phase_1}, system \eqref{averaged}) (black color)
and physical model (see Fig.~\ref{physical_sim_model}, system \eqref{prev diff eq}) (red color)
after transient processes have different phases in the locked states.

Here VCO free-running frequency $\omega_2^{free}=2*\pi*400000-600000$;
initial states of filters are all zero: $x(0) = x_1(0)\equiv x_2{0} = 0$
(i.e. $\alpha_0(t,x(0)) = \alpha_{1,2}(t,x_{1,2}(0)) \equiv 0$).

\begin{figure}
  \centering\includegraphics[width=0.23\textwidth]{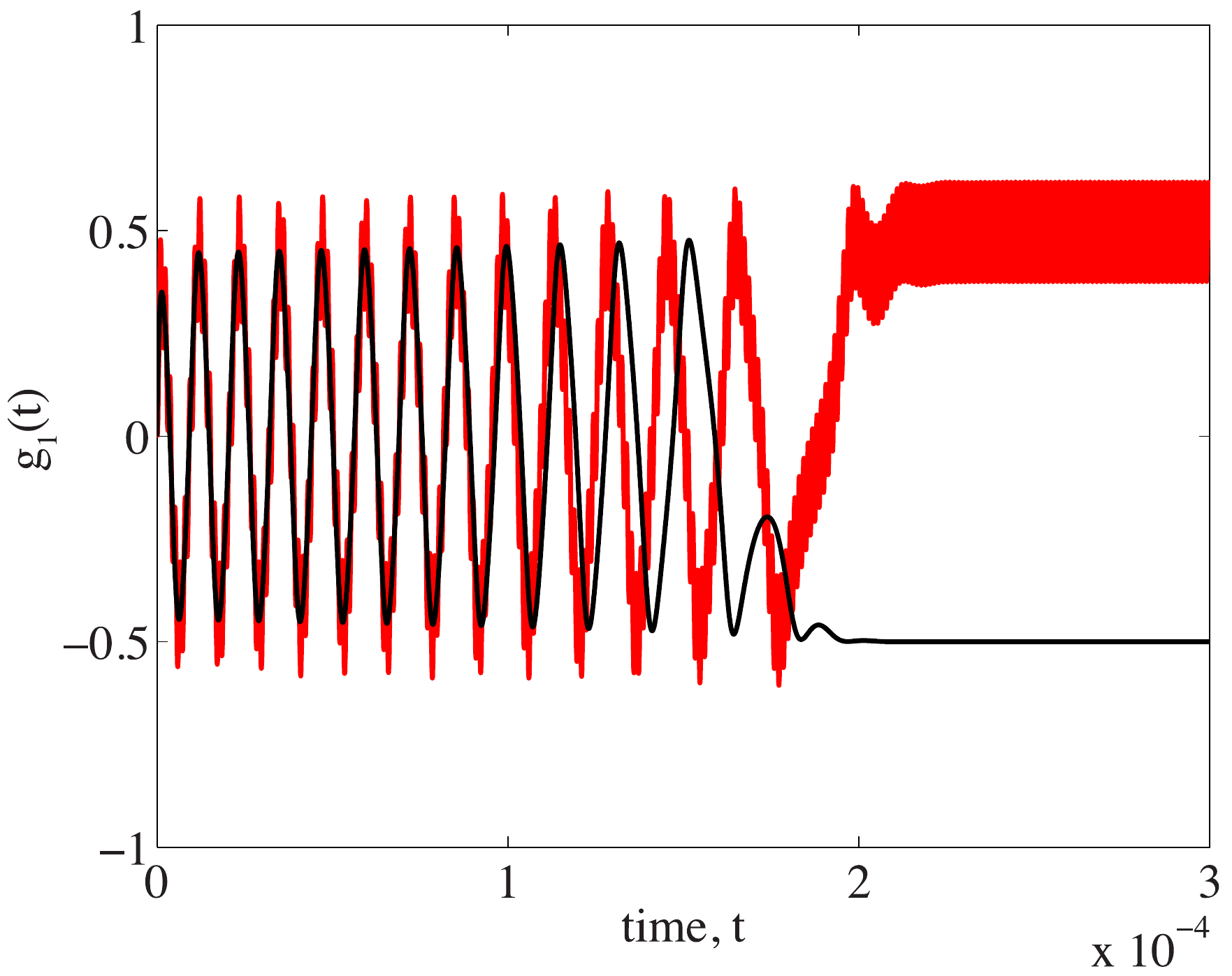}
  \centering\includegraphics[width=0.23\textwidth]{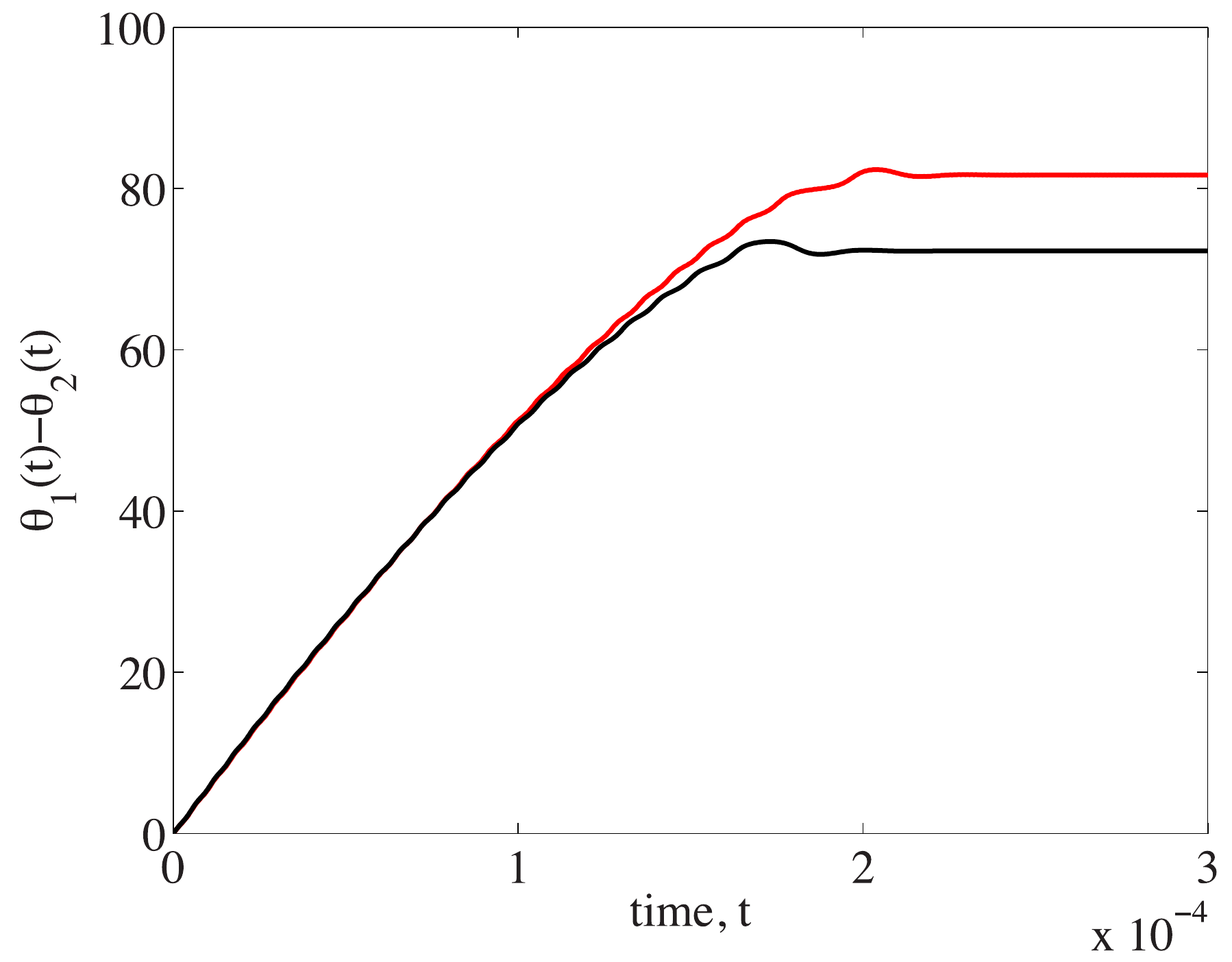}
  \caption{Low-pass filter outputs $g_1(t)$ and phase difference $\theta_\Delta(t)$ for
  averaged model \eqref{averaged} (black)
  and physical model (red) in Fig.~\ref{costas_out}.}
  \label{averaged_non_averaged}
\end{figure}
\end{example}

\begin{example}[initial states of the low-pass filters]
In~Fig.~\ref{lpf_ics}
it is shown that Assumption~\ref{as-lpf-initstate} may not be valid:
while physical model (see Fig.~\ref{physical_sim_model}, system \eqref{prev diff eq})
with zero initial states of low-pass filters acquires lock (black),
the same physical model with nonzero initial states of low-pass filters
is out of lock (red).

Here VCO free-running frequency $\omega_2^{free}= 2*\pi*400000-2$,
initial state of the loop filter is zero: $x(0)=0$ (i.e. $\alpha_0(t,x(0)) \equiv 0$),
Initial states of the low-pass filters are $x_1(0)=x_2(0)=0$ (black)
and $x_1(0)=0.02$, $x_2(0)=0$ (red).

\begin{figure}
  \centering
  \includegraphics[width=0.23\textwidth]{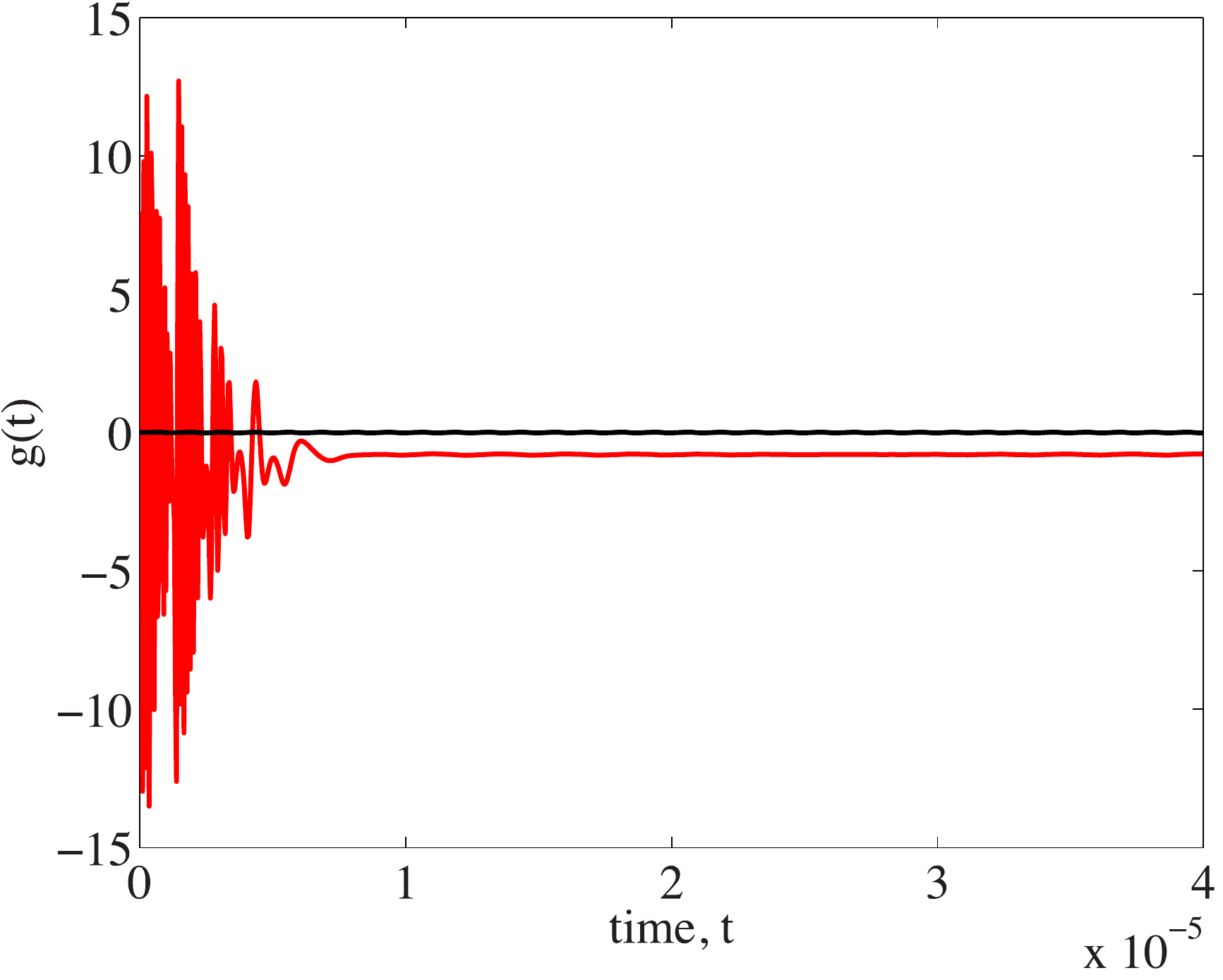}
  \includegraphics[width=0.23\textwidth]{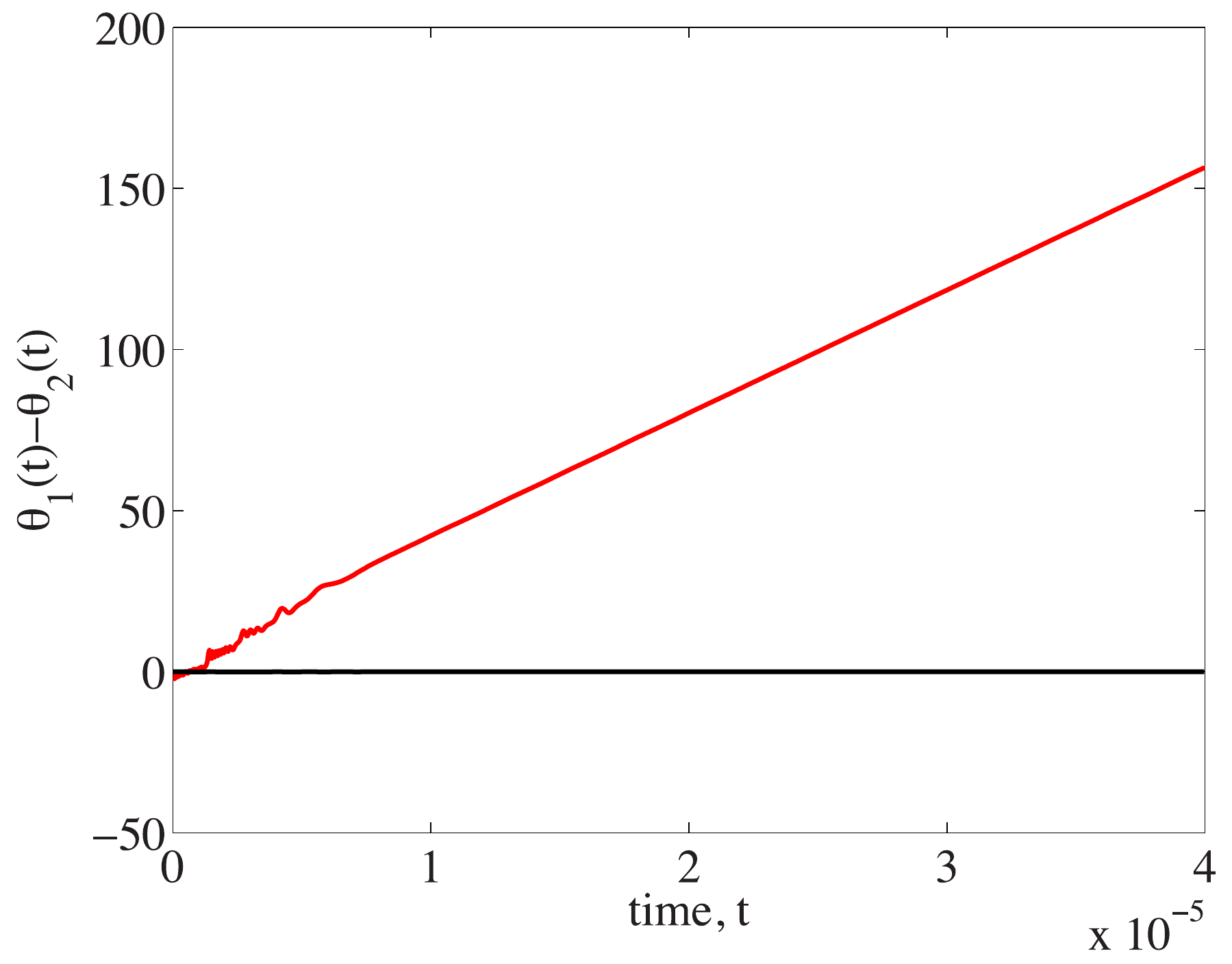}
  \caption{
  Loop filter output $g(t)$ for
  physical model (black) with zero initial states of low-pass filters,
  physical model (red) with nonzero initial states of low-pass filters.}
  \label{lpf_ics}
\end{figure}
\end{example}

\begin{example}[data signal]
In~Fig.~\ref{data_signal}
it is shown that Assumption~\ref{as-lpf-data} may not be valid:
while the simplified mathematical model in the signal space
(see Fig.~\ref{physical_sim_model} with $m(t)\equiv 1$, system \eqref{diff_eq_0})
acquires lock (black),
physical model with (see Fig.~\ref{physical_sim_model}, system \eqref{prev diff eq})
with periodic data is out of lock (red).

Here VCO free-running frequency $\omega_2^{free}=3.2*10^6$,
initial states of filters are all zero: $x(0)=x_1(0)=x_2(0)=0$,
data signal is periodic: $m(t)=\sign\sin(10^5*2*\pi*t)$.
\begin{figure}
  \centering
  \includegraphics[width=0.23\textwidth]{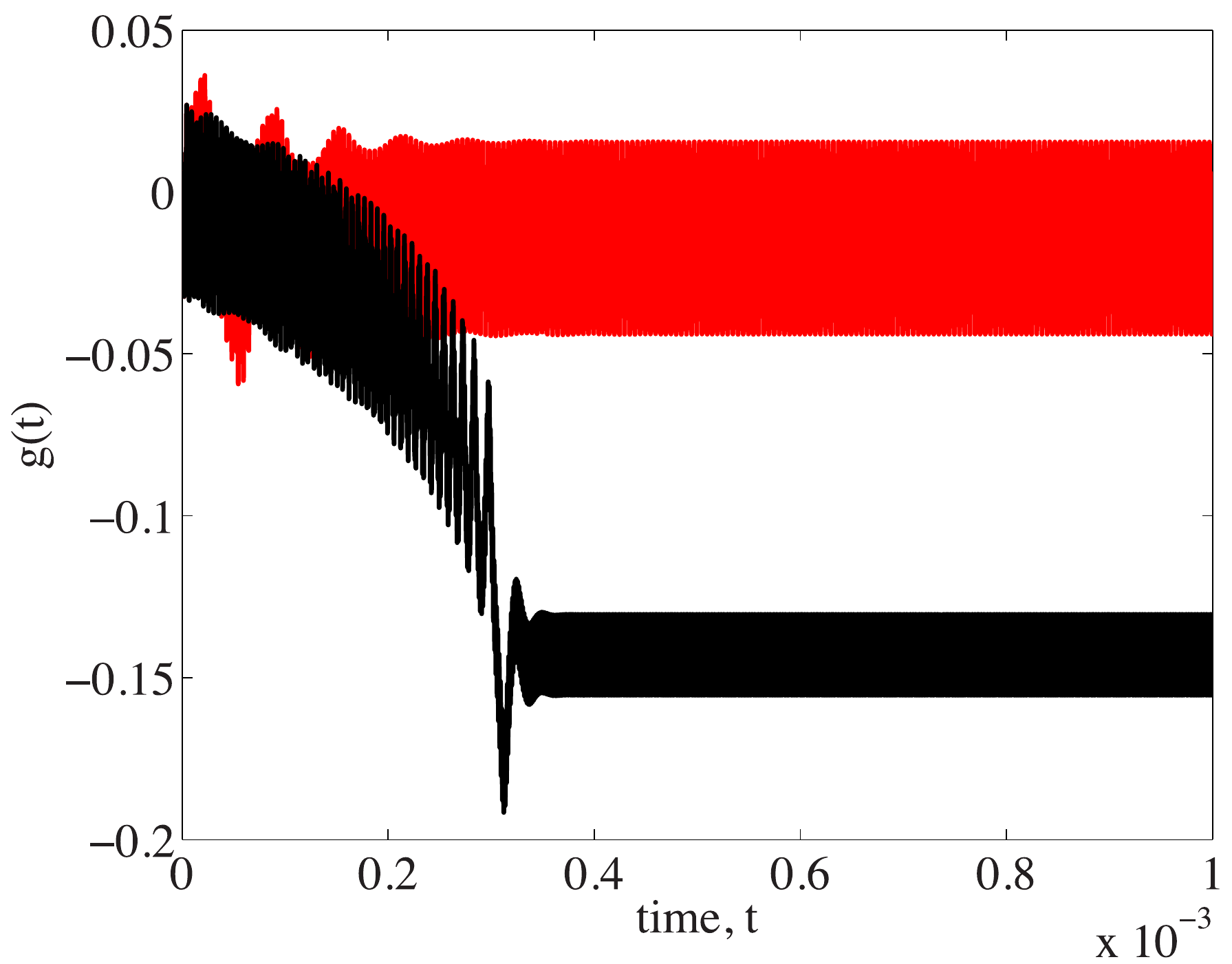}
  \includegraphics[width=0.23\textwidth]{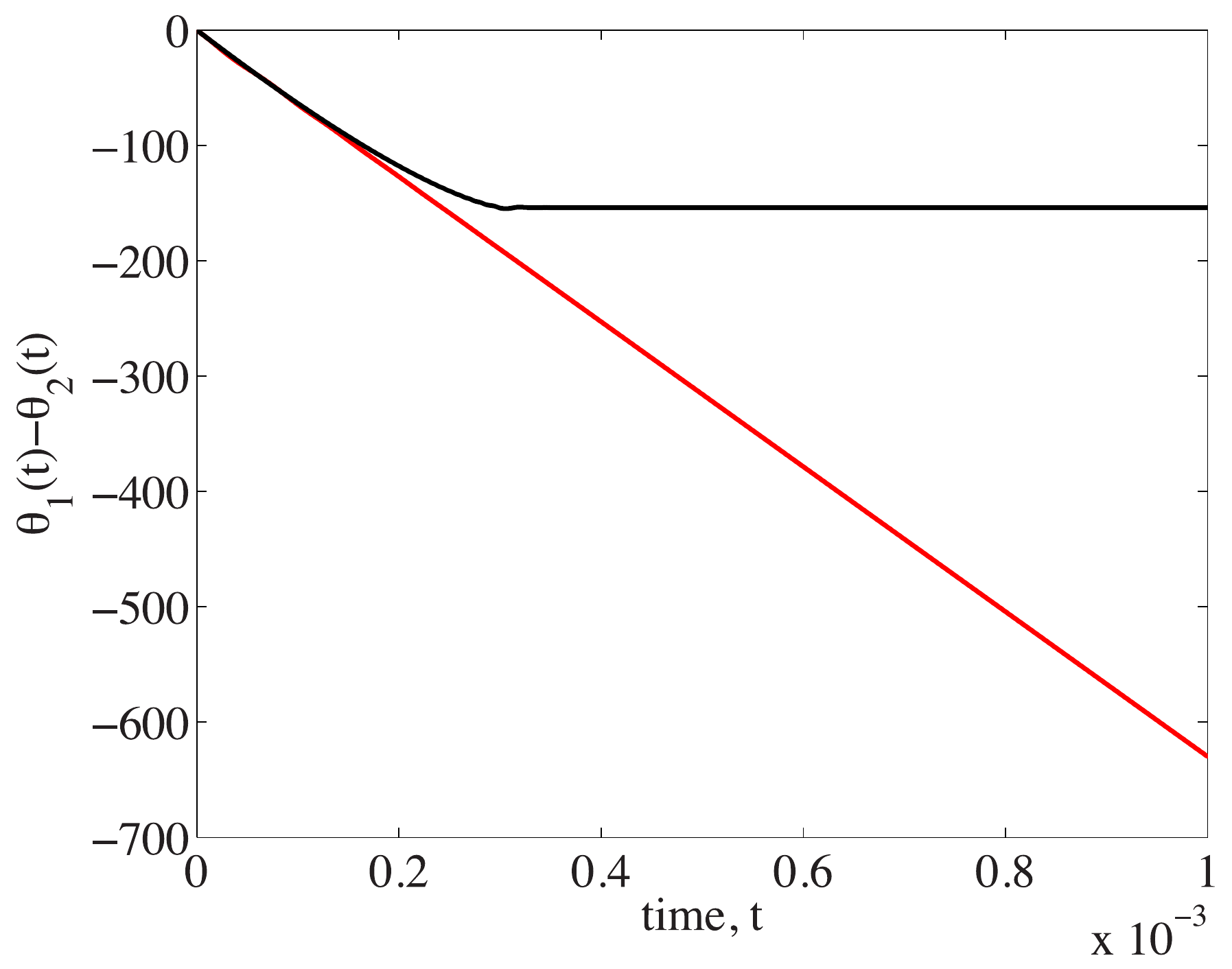}
  \caption{
  Loop filter output $g(t)$ for
  physical model (red) with periodic data signal,
  physical model (black) with constant data signal $m(t)\equiv 1$.}
  \label{data_signal}
\end{figure}
\end{example}

\begin{example}[low-pass filters]
In~Fig.~\ref{arm_filters_instability}
it is shown that the consideration of the ideal low-pass filters
(Assumptions~\ref{as-twice-frequency},\ref{as-lpf-initstate}, and \ref{as-lpf-data})
may lead to wrong conclusions:
it is shown that low-pass filters may affect the stability of models in signal's phase space.
While simplified mathematical model in signal's phase space
(see Fig.~\ref{costas_phase_1}, system \eqref{averaged}) (black)
and simplified mathematical model in the signal space
(see Fig.~\ref{physical_sim_model} with $m(t) \equiv 1$, system \eqref{diff_eq_0}) (green)
are out of lock,
the classic mathematical model in signal's phase space
(see Fig.~\ref{costas_phase_1_short}, system \eqref{final_system}), where low-pass filters are not taken into account,
acquires lock (red).

Here VCO free-running frequency $\omega_2^{free}=2*\pi*400000-500000$,
initial states of filters are zero: $x(0)=x_1(0)=x_2(0)=0$,
no data is being transmitted: $m(t) \equiv 1$.

\begin{figure}
  \includegraphics[width=0.23\textwidth]{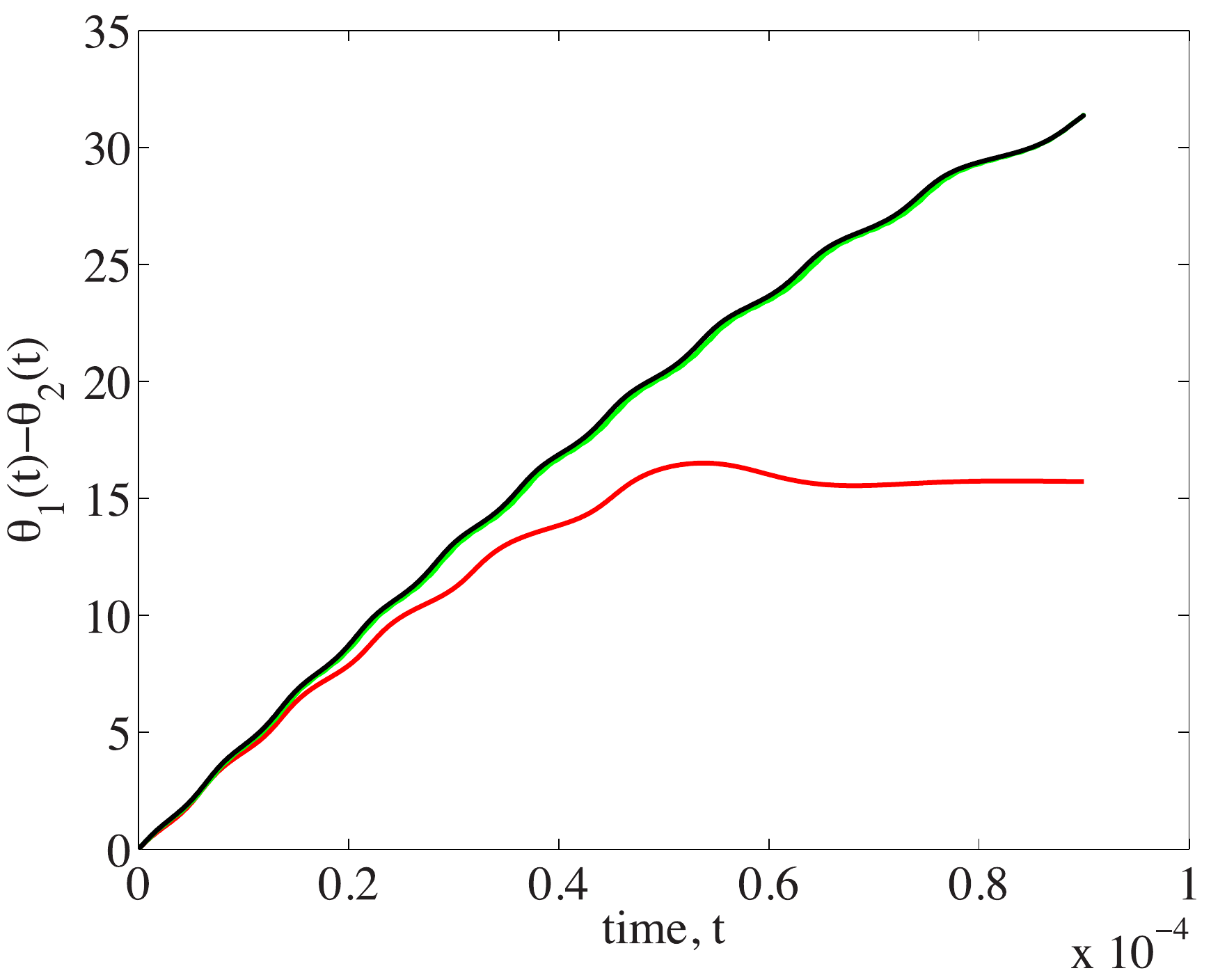}
  \includegraphics[width=0.23\textwidth]{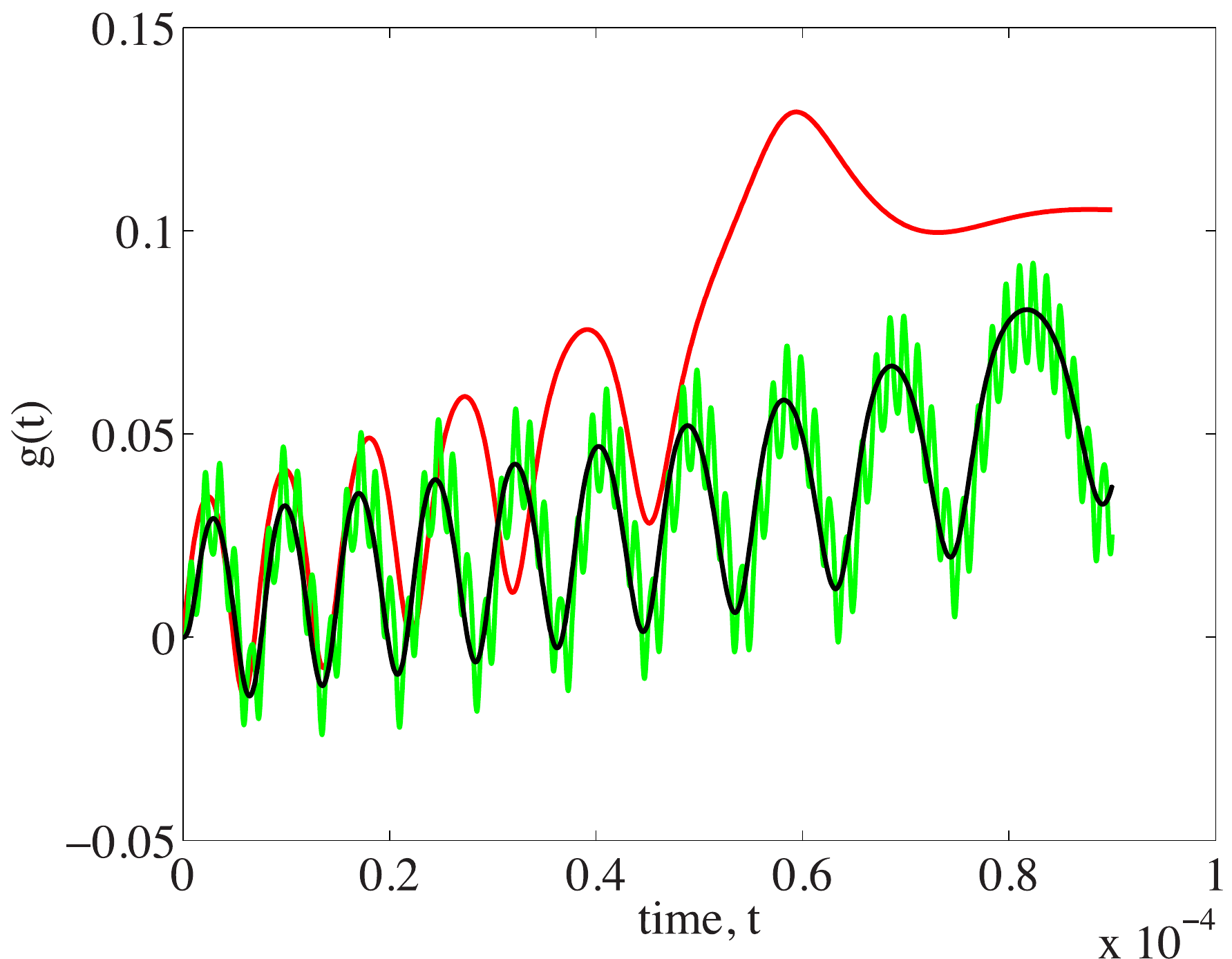}
  \caption{Loop filter output $g(t)$ and phase difference $\theta_\Delta(t)$ for
  signal's phase space model (black) with low-pass filters,
  classic signal's phase space model (red) without low-pass filters,
  and signal space model (green).}
  \label{arm_filters_instability}
\end{figure}
\end{example}

\begin{example}[initial state of the loop filter]
In~Fig.~\ref{loop_filter_ics}
it is shown that Assumption~\ref{as-lf-initstate} may not be valid:
while the physical model (see Fig.~\ref{physical_sim_model}, system \eqref{prev diff eq})
with zero output of loop filter acquires lock (black),
the same physical model with nonzero output of loop filter is out of lock (red).

Here initial states of low-pass filters are zero: $x_1(0)=x_2(0)=0$;
VCO free-running frequency $\omega_2^{free}=2*\pi*400000-10$;
initial loop filter state is $x(0)=-0.00001$ (red)
and $x(0)=0$ (black).

\begin{figure}
  \centering\includegraphics[width=0.23\textwidth]{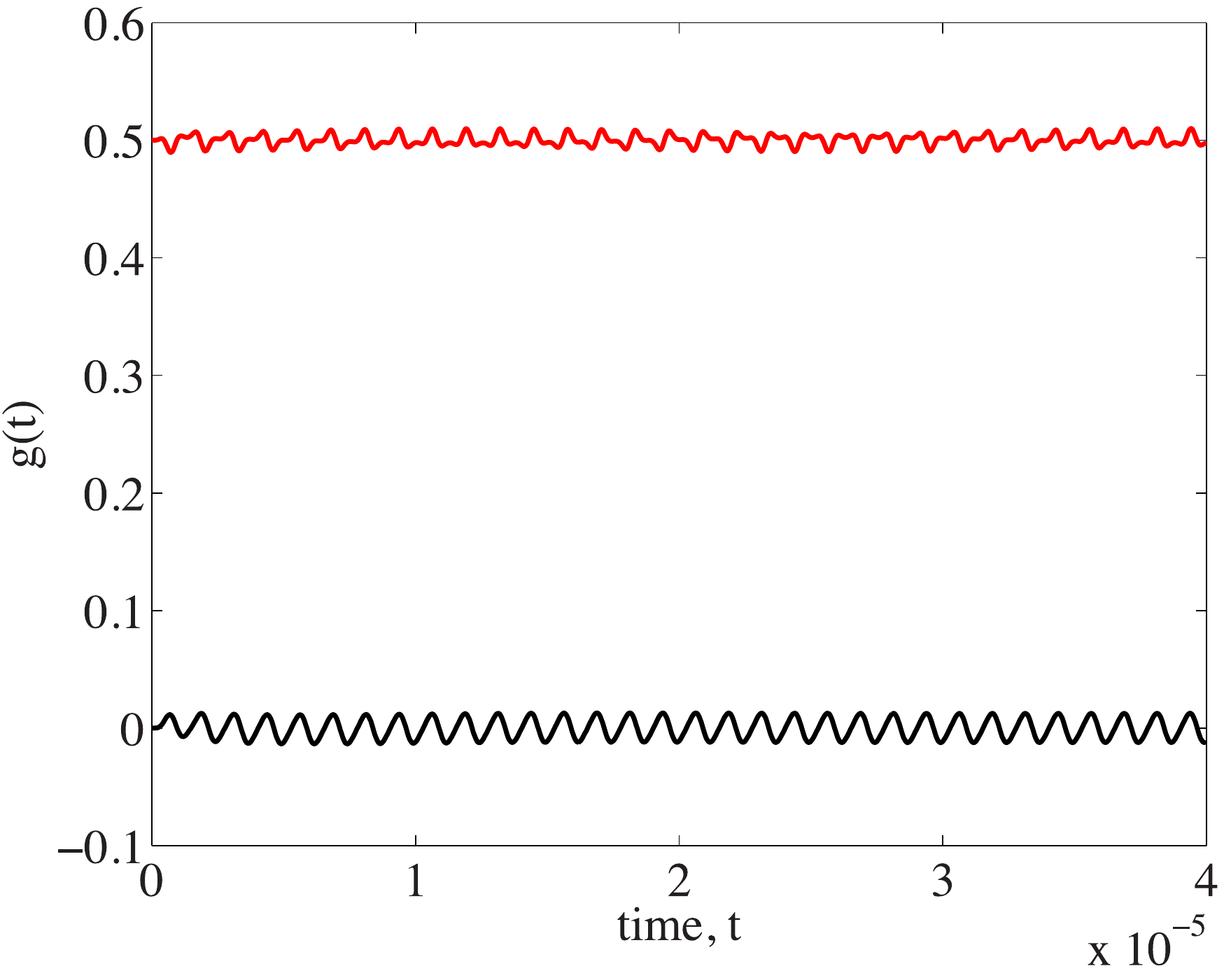}
  \centering\includegraphics[width=0.23\textwidth]{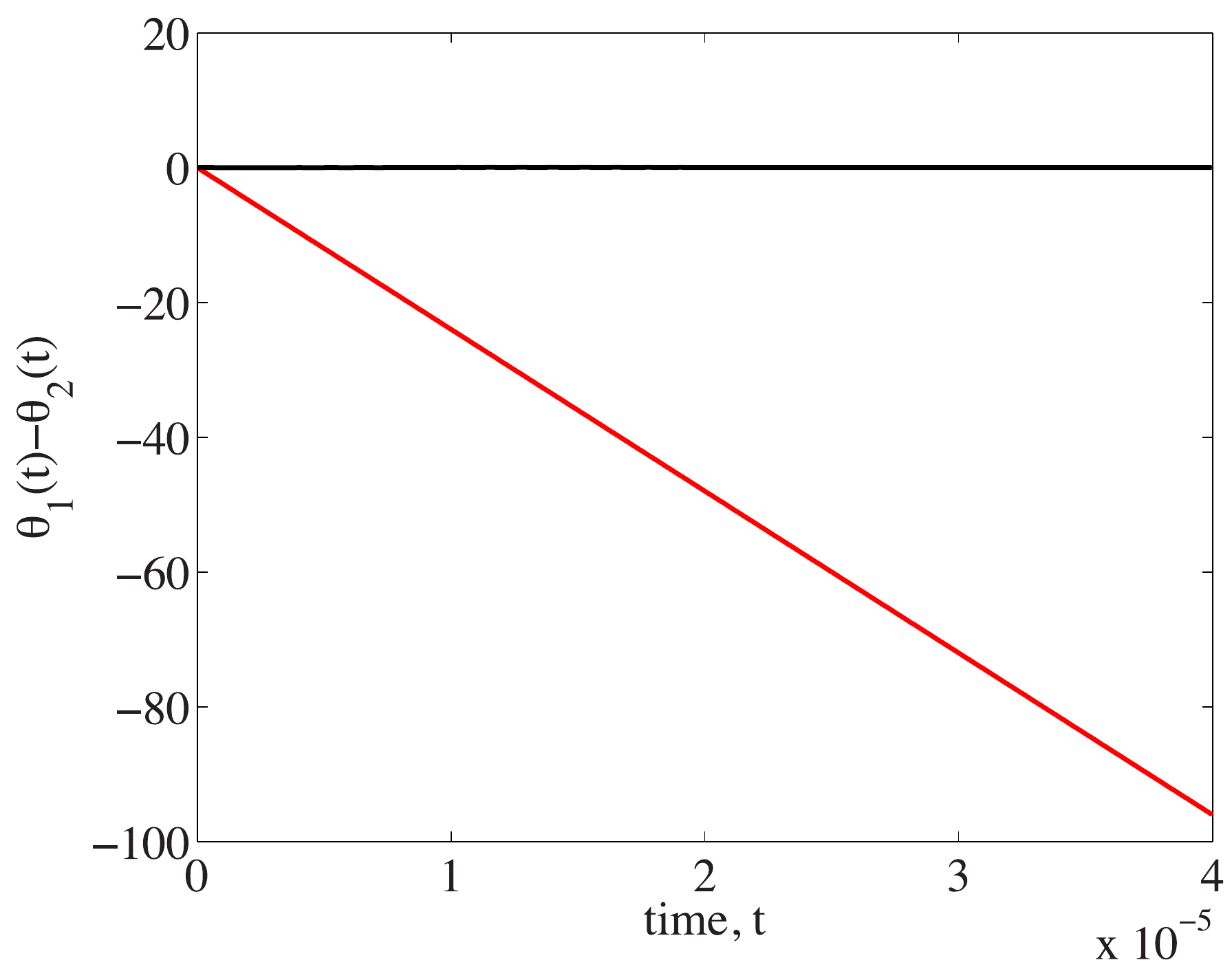}
  \caption{
  Loop filter output $g(t)$ and phase difference $\theta_\Delta(t)$ for
  physical model with zero initial state of loop filter (black),
  physical model with nonzero initial state of loop filter (red).
 }
  \label{loop_filter_ics}
\end{figure}
\end{example}

\begin{example}[numerical integration parameters]
In~Fig.~\ref{semistable}
it is shown that standard simulation of the loop may not be valid: while
the classic mathematical model in signal's phase space
(Fig.~\ref{costas_phase_1_short} or system \eqref{final_system}),
simulated in Simulink with predefined integration parameters: 'max step size' set to '1e-3',
is out of lock  (black),
the same model simulated in Simulink
with default integration parameters: 'max step size' set to 'auto', acquires lock (red).
Here Matlab chooses step from $5 \cdot 10^{-3}$ to $9\cdot10^{-2}$;
for fixed step $2\cdot 10^{-2}$ the model acquires lock,
for fixed step $1\cdot 10^{-2}$ the model doesn't acquire lock.

Here the initial loop filter state output is $x(0) = 0.0125$;
VCO free-running frequency $\omega_2^{free}=10000-89.45$;
VCO input gain $L=1000$;
initial phase shift $\theta_\Delta(0)=-3.4035$.

\begin{figure}
  \centering\includegraphics[width=0.23\textwidth]{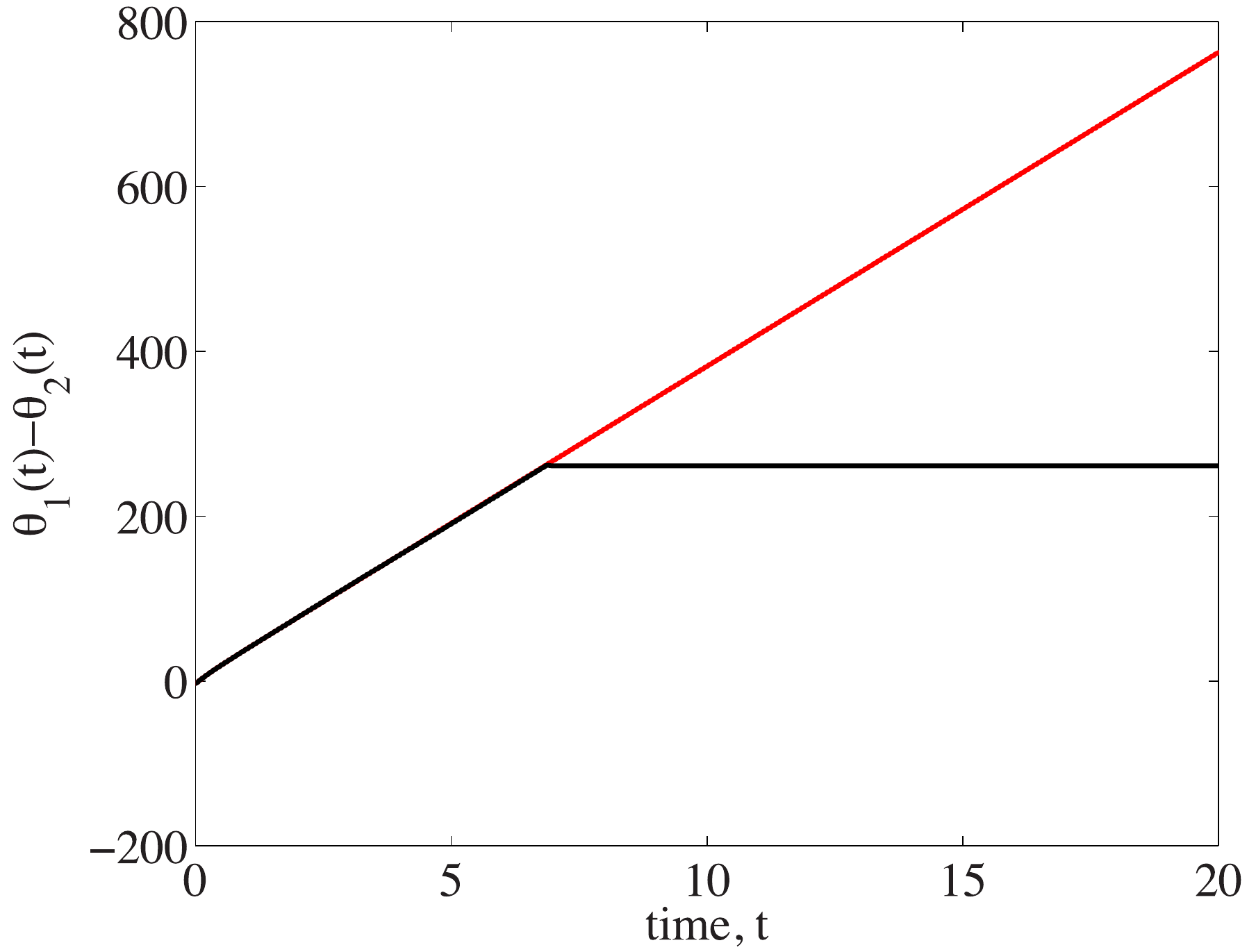}
  \centering\includegraphics[width=0.23\textwidth]{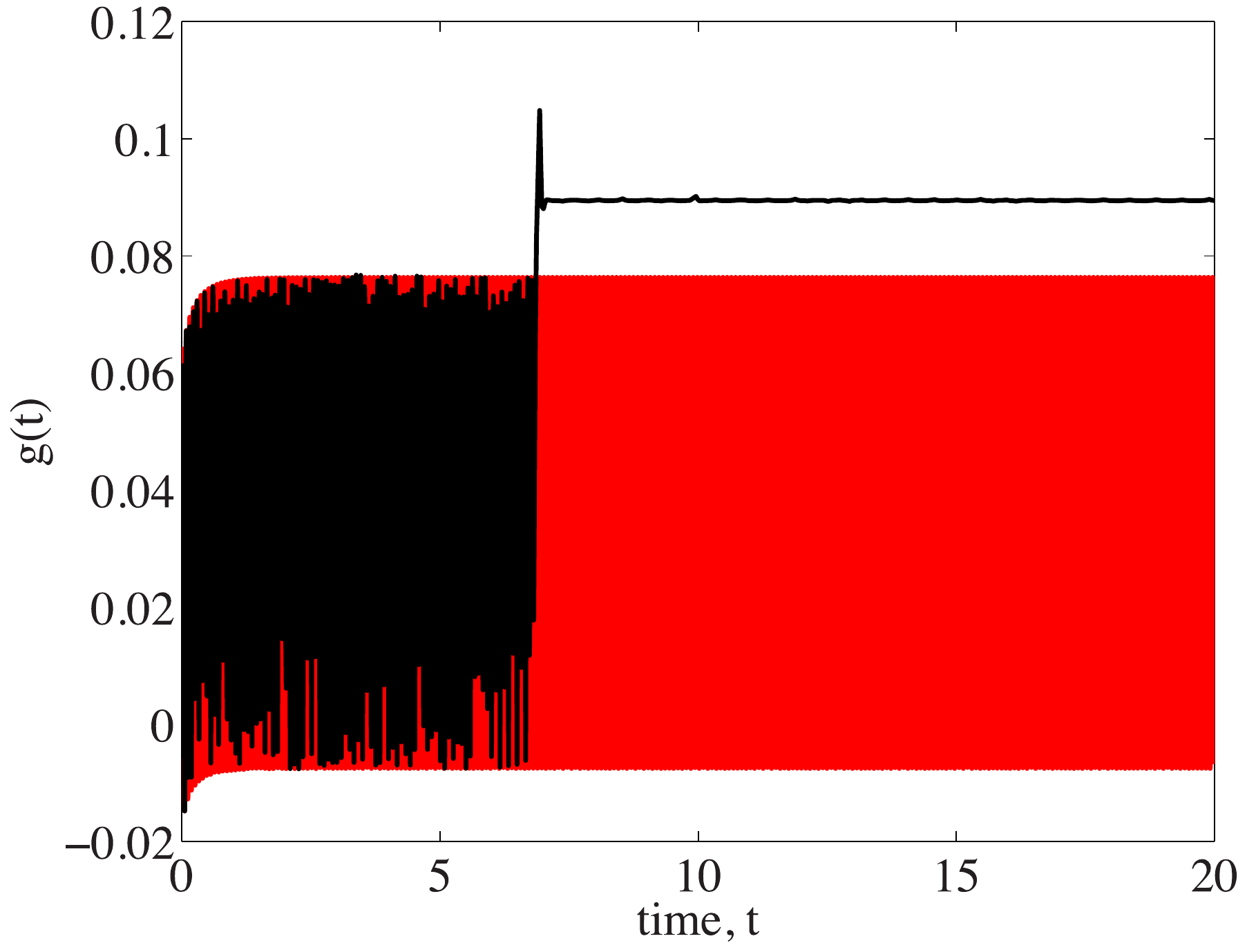}
  \caption{Filter outputs: default integration parameters in Simulink
  'max step size' set to 'auto' (black curve);
  Parameters configured manually 'max step size' set to '1e-3' (red curve).}
  \label{semistable}
\end{figure}
\end{example}

\begin{figure}
  \centering\includegraphics[width=0.30\textwidth]{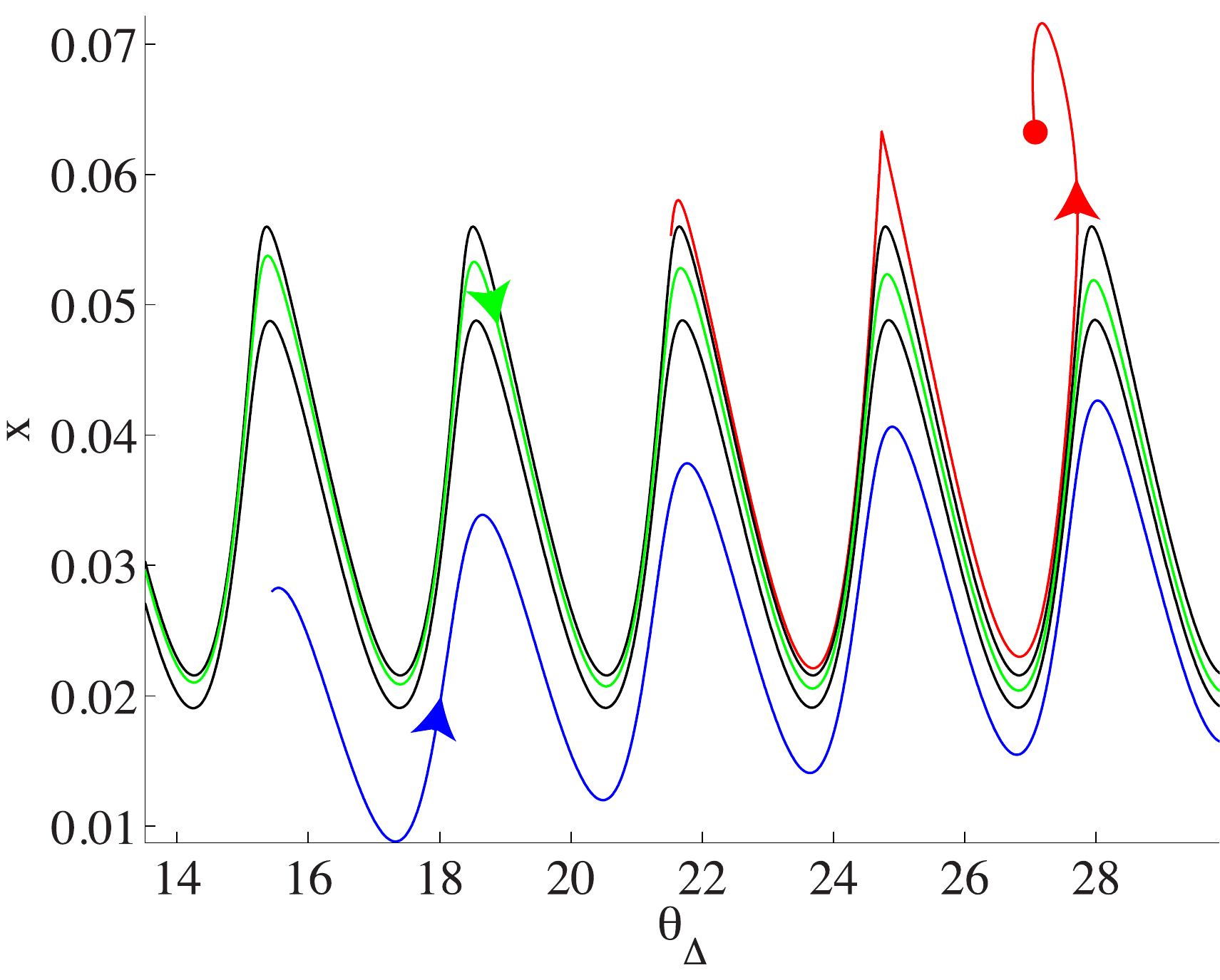}
  \caption{Phase portrait: coexistence of stable and unstable periodic solutions.}
  \label{semistable_phase}
\end{figure}

Consider now the corresponding phase portrait (see Fig.~\ref{semistable_phase}).
Here the red trajectory tends to a stable equilibrium (red dot).
Lower and higher black trajectories are stable and unstable limit cycles.
The blue trajectory tends to a stable periodic trajectory (lower black periodic curve) and
in this case the model does not acquire lock.
All trajectories between black trajectories (see green trajectory) tend to the stable lower black trajectory.

If the gap between stable and unstable trajectories (black lines) is smaller than the discretization step,
the numerical procedure may slip through the stable trajectory
(blue trajectory may step over the black and green lines and begins to be attracted to the red dot).
In other words, the simulation may show that the Costas loop acquires lock
although in reality it does not.
The considered case corresponds to the coexisting attractors (one of which is a hidden oscillation)
and the bifurcation of birth of a semistable trajectory \cite{LeonovK-2013-IJBC}.

Note, that only trajectories (red) above the unstable limit cycle is attracted to the equilibrium.
Hence $\omega_\Delta = 89.45$ does not belong to the pull-in range.

\section{Conclusion}
In this paper a short survey on derivation
of nonlinear mathematical models of BPSK Costas loop, used for pre-design and post-design analysis,
is presented. It has been shown that
1) the consideration of simplified mathematical models,
constructed intuitively,
and
2) the application of non rigorous methods of analysis (e.g., simulation and linearization)
can lead to wrong conclusions concerning the operability of the Costas loop \emph{physical model}.
Similar result can be obtained for the classic PLL and QPSK Costas loop
(see, e.g., \cite{KuznetsovLYY-2014-IFAC,KuznetsovKLNYY-2015-ISCAS,KuznetsovKLNYY-2014-ICUMT-QPSK}).

While the Costas loop is a nonlinear control system and for its analysis
it is essential to apply various stability criteria developed in the control theory,
their direct application to the PLL-based models is impossible,
because such criteria usually are not adapted for the cylindrical phase space.
In the tutorial {\it Phase Locked Loops: a  Control Centric Tutorial},
presented at the {\it American Control Conference 2002},
D.~Abramovitch wrote that
``{\it The general theory of PLLs and ideas on how to make them even more useful seems
to cross into the controls literature only rarely}'' \cite{Abramovitch-2002}.

Corresponding modifications of the classic stability criteria
for the rigorous analytical analysis of nonlinear PLL-based model in the cylindrical phase space
had been well developed in 197x-199x in \cite{GeligLY-1978,LeonovRS-1992,LeonovBSh-1996}. %,LeonovPS-1996
See also some recent works on nonlinear methods for the analysis of PLL-based models
\cite{SuarezQ-2003,YakubovichLG-2004,Margaris-2004,Leonov-2006-ARC,KudrewiczW-2007,ChiconeH-2013,LeonovKS-2009,LeonovK-2014}.
One of the reason why these works have been remained almost unnoticed by the engineers
may be that they are mostly written in
the language of the control theory and the theory of dynamical systems,
and, thus, may not be well adapted to the terms and objects used in the engineering practice of phase-locked loops.
Another possible reason is that
\,``{\it nonlinear analysis techniques are well beyond the scope of most undergraduate courses in communication theory}''
\cite{Tranter-2010-book}. %[p.~1]

The examples, considered in the paper, are the motivation to apply rigorous analytical methods
for the analysis of PLL-based loop nonlinear models.

\section*{\uppercase{Acknowledgements}}
This work was supported by
Saint-Petersburg State University (project 6.39.416.2014),
the Leading Scientific Schools programm (3384.2014.1),
Russian Scientific Foundation (14-21-00041, sec.II-III).

\bibliographystyle{IEEEtran}
%\bibliography{C:/Dropbox/bib/bib_nk,C:/Dropbox/bib/bib_leonov,C:/Dropbox/bib/bib_full,C:/Dropbox/bib/bib_pll,C:/Dropbox/bib/bib-gly}

\begin{thebibliography}{10}
\providecommand{\url}[1]{#1}
\csname url@rmstyle\endcsname
\providecommand{\newblock}{\relax}
\providecommand{\bibinfo}[2]{#2}
\providecommand\BIBentrySTDinterwordspacing{\spaceskip=0pt\relax}
\providecommand\BIBentryALTinterwordstretchfactor{4}
\providecommand\BIBentryALTinterwordspacing{\spaceskip=\fontdimen2\font plus
\BIBentryALTinterwordstretchfactor\fontdimen3\font minus
  \fontdimen4\font\relax}
\providecommand\BIBforeignlanguage[2]{{%
\expandafter\ifx\csname l@#1\endcsname\relax
\typeout{** WARNING: IEEEtran.bst: No hyphenation pattern has been}%
\typeout{** loaded for the language `#1'. Using the pattern for}%
\typeout{** the default language instead.}%
\else
\language=\csname l@#1\endcsname
\fi
#2}}

\bibitem{Costas-1956}
J.~{C}ostas, ``Synchoronous communications,'' in \emph{Proc. IRE}, vol.~44,
  1956, pp. 1713--1718.

\bibitem{Costas-1962-patent}
J.~P. Costas, ``Receiver for communication system,'' July 1962, {US} Patent
  3,047,659.

\bibitem{Best-2007}
R.~Best, \emph{Phase-Lock Loops: Design, Simulation and Application},
  6th~ed.\hskip 1em plus 0.5em minus 0.4em\relax McGraw-Hill, 2007.

\bibitem{KiharaOE-2002}
M.~Kihara, S.~Ono, and P.~Eskelinen, \emph{Digital Clocks for Synchronization
  and Communications}.\hskip 1em plus 0.5em minus 0.4em\relax Artech House,
  2002.

\bibitem{Bianchi-2005-book}
G.~Bianchi, \emph{Phase-Locked Loop Synthesizer Simulation}.\hskip 1em plus
  0.5em minus 0.4em\relax McGraw-Hill, 2005.

\bibitem{PedersonM-2008-book}
D.~Pederson and K.~Mayaram, \emph{Analog Integrated Circuits for Communication:
  Principles, Simulation and Design}.\hskip 1em plus 0.5em minus 0.4em\relax
  Springer, 2008.

\bibitem{Tranter-2010-book}
W.~Tranter, T.~Bose, and R.~Thamvichai, \emph{Basic Simulation Models of Phase
  Tracking Devices Using {MATLAB}}, ser. Synthesis lectures on
  communications.\hskip 1em plus 0.5em minus 0.4em\relax Morgan \& Claypool,
  2010.

\bibitem{Talbot-2012-book}
D.~Talbot, \emph{Frequency Acquisition Techniques for Phase Locked
  Loops}.\hskip 1em plus 0.5em minus 0.4em\relax Wiley, 2012.

\bibitem{Viterbi-1966}
A.~Viterbi, \emph{Principles of coherent communications}.\hskip 1em plus 0.5em
  minus 0.4em\relax New York: McGraw-Hill, 1966.

\bibitem{PiqueiraM-2003}
J.~Piqueira and L.~Monteiro, ``Considering second-harmonic terms in the
  operation of the phase detector for second-order phase-locked loop,''
  \emph{IEEE Transactions On Circuits And Systems-I}, vol.~50, no.~6, pp.
  805--809, 2003.

\bibitem{KuznetsovKLNYY-2015-ISCAS}
N.~Kuznetsov, O.~Kuznetsova, G.~Leonov, P.~Neittaanmaki, M.~Yuldashev, and
  R.~Yuldashev, ``Limitations of the classical phase-locked loop analysis,'' in
  \emph{International Symposium on Circuits and Systems (ISCAS)}.\hskip 1em
  plus 0.5em minus 0.4em\relax IEEE, 2015, accepted.

\bibitem{MitropolskyB-1961}
Y.~Mitropolsky and N.~Bogolubov, \emph{Asymptotic Methods in the Theory of
  Non-Linear Oscillations}.\hskip 1em plus 0.5em minus 0.4em\relax New York:
  Gordon and Breach, 1961.

\bibitem{LeonovKYY-2012-TCASII}
G.~A. Leonov, N.~V. Kuznetsov, M.~V. Yuldahsev, and R.~V. Yuldashev,
  ``Analytical method for computation of phase-detector characteristic,''
  \emph{IEEE Transactions on Circuits and Systems - II: Express Briefs},
  vol.~59, no.~10, pp. 633--647, 2012.

\bibitem{LeonovKYY-2015-SIGPRO}
G.~A. Leonov, N.~V. Kuznetsov, M.~V. Yuldashev, and R.~V. Yuldashev,
  ``Nonlinear dynamical model of {C}ostas loop and an approach to the analysis
  of its stability in the large,'' \emph{Signal processing}, vol. 108, pp.
  124--135, 2015.

\bibitem{Olson-1975}
M.~Olson, ``False-lock detection in {C}ostas demodulators,'' \emph{Aerospace
  and Electronic Systems, IEEE Transactions on}, vol. AES-11, no.~2, pp.
  180--182, 1975.

\bibitem{Stensby-1989}
J.~Stensby, ``False lock and bifurcation in {C}ostas loops,'' \emph{SIAM
  Journal on Applied Mathematics}, vol.~49, no.~2, pp. pp. 420--431, 1989.

\bibitem{LeonovK-2014}
G.~A. Leonov and N.~V. Kuznetsov, \emph{Nonlinear Mathematical Models Of
  Phase-Locked Loops. Stability and Oscillations}.\hskip 1em plus 0.5em minus
  0.4em\relax Cambridge Scientific Press, 2014.

\bibitem{AscheidM-1982}
G.~Ascheid and H.~Meyr, ``Cycle slips in phase-locked loops: A tutorial
  survey,'' \emph{Communications, IEEE Transactions on}, vol.~30, no.~10, pp.
  2228--2241, 1982.

\bibitem{ErshovaL-1983}
O.~B. Ershova and G.~A. Leonov, ``Frequency estimates of the number of cycle
  slidings in phase control systems,'' \emph{Avtomat. Remove Control}, vol.~44,
  no.~5, pp. 600--607, 1983.

\bibitem{Goyal-2006}
P.~Goyal, X.~Lai, and J.~Roychowdhury, ``A fast methodology for
  first-time-correct design of {PLL}s using nonlinear phase-domain {VCO}
  macromodels,'' in \emph{Proceedings of the 2006 Asia and South Pacific Design
  Automation Conference}, 2006, pp. 291--296.

\bibitem{KuznetsovLYY-2011-IEEE}
N.~V. Kuznetsov, G.~A. Leonov, M.~V. Yuldashev, and R.~V. Yuldashev,
  ``Analytical methods for computation of phase-detector characteristics and
  {PLL} design,'' in \emph{ISSCS 2011 - International Symposium on Signals,
  Circuits and Systems, Proceedings}.\hskip 1em plus 0.5em minus 0.4em\relax
  IEEE, 2011, pp. 7--10.

\bibitem{ourPatent-2013-fin}
N.~V. Kuznetsov, G.~A. Leonov, P.~Neittaanmaki, M.~V. Yuldashev, and R.~V.
  Yuldashev, ``Method and system for modeling {C}ostas loop feedback for fast
  mixed signals,'' Finland Patent, 2013, {P}atent application FI20130124.

\bibitem{Abramovitch-2002}
D.~Abramovitch, ``Phase-locked loops: A control centric tutorial,'' in
  \emph{American Control Conf. Proc.}, vol.~1, 2002, pp. 1--15.

\bibitem{LeonovK-2013-IJBC}
G.~A. Leonov and N.~V. Kuznetsov, ``Hidden attractors in dynamical systems.
  {F}rom hidden oscillations in {H}ilbert-{K}olmogorov, {A}izerman, and
  {K}alman problems to hidden chaotic attractors in {C}hua circuits,''
  \emph{International Journal of Bifurcation and Chaos}, vol.~23, no.~1, 2013,
  {a}rt. no. 1330002.

\bibitem{KuznetsovLYY-2014-IFAC}
N.~Kuznetsov, G.~Leonov, M.~Yuldashev, and R.~Yuldashev, ``Nonlinear analysis
  of classical phase-locked loops in signal's phase space,'' \emph{IFAC
  Proceedings Volumes (IFAC-PapersOnline)}, vol.~19, no.~1, pp. 8253--8258,
  2014.

\bibitem{KuznetsovKLNYY-2014-ICUMT-QPSK}
N.~Kuznetsov, O.~Kuznetsova, G.~Leonov, P.~Neittaanmaki, M.~Yuldashev, and
  R.~Yuldashev, ``Simulation of nonlinear models of {QPSK} {C}ostas loop in
  {M}atlab {S}imulink,'' in \emph{Ultra Modern Telecommunications and Control
  Systems and Workshops (ICUMT), 2014 6th International Congress on}.\hskip 1em
  plus 0.5em minus 0.4em\relax IEEE, 2014, pp. 66--71.

\bibitem{GeligLY-1978}
A.~Gelig, G.~Leonov, and V.~Yakubovich, \emph{Stability of Nonlinear Systems
  with Nonunique Equilibrium (in Russian)}.\hskip 1em plus 0.5em minus
  0.4em\relax Nauka, 1978.

\bibitem{LeonovRS-1992}
G.~A. Leonov, V.~Reitmann, and V.~B. Smirnova, \emph{Nonlocal Methods for
  Pendulum-like Feedback Systems}.\hskip 1em plus 0.5em minus 0.4em\relax
  Stuttgart-Leipzig: Teubner Verlagsgessel\-schaft, 1992.

\bibitem{LeonovBSh-1996}
G.~A. Leonov, I.~M. Burkin, and A.~I. Shepelyavy, \emph{Frequency Methods in
  Oscillation Theory}.\hskip 1em plus 0.5em minus 0.4em\relax Dordretch:
  Kluwer, 1996.

\bibitem{SuarezQ-2003}
A.~Suarez and R.~Quere, \emph{Stability Analysis of Nonlinear Microwave
  Circuits}.\hskip 1em plus 0.5em minus 0.4em\relax Artech House, 2003.

\bibitem{YakubovichLG-2004}
V.~A. Yakubovich, G.~A. Leonov, and A.~K. Gelig, \emph{Stability of Stationary
  Sets in Control Systems with Discontinuous Nonlinearities}.\hskip 1em plus
  0.5em minus 0.4em\relax Singapure: World Scientific, 2004.

\bibitem{Margaris-2004}
N.~Margaris, \emph{Theory of the Non-Linear Analog Phase Locked Loop}.\hskip
  1em plus 0.5em minus 0.4em\relax New Jersey: Springer Verlag, 2004.

\bibitem{Leonov-2006-ARC}
G.~A. Leonov, ``Phase-locked loops. {T}heory and application,''
  \emph{Automation and Remote Control}, vol.~10, pp. 47--55, 2006.

\bibitem{KudrewiczW-2007}
J.~Kudrewicz and S.~Wasowicz, \emph{Equations of phase-locked loop. Dynamics on
  circle, torus and cylinder}.\hskip 1em plus 0.5em minus 0.4em\relax World
  Scientific, 2007.

\bibitem{ChiconeH-2013}
C.~Chicone and M.~Heitzman, ``Phase-locked loops, demodulation, and averaging
  approximation time-scale extensions,'' \emph{SIAM J. Applied Dynamical
  Systems}, vol.~12, no.~2, pp. 674--721, 2013.

\bibitem{LeonovKS-2009}
G.~A. Leonov, N.~V. Kuznetsov, and S.~M. Seledzhi, \emph{Automation control -
  Theory and Practice}.\hskip 1em plus 0.5em minus 0.4em\relax In-Tech, 2009,
  ch. Nonlinear Analysis and Design of Phase-Locked Loops, pp. 89--114.

\end{thebibliography}

\newcommand{\noopsort}[1]{} \newcommand{\printfirst}[2]{#1}
  \newcommand{\singleletter}[1]{#1} \newcommand{\switchargs}[2]{#2#1}

\end{document}